\documentclass[11pt]{amsart}
\usepackage{amsxtra}
\usepackage{amssymb}
\usepackage[mathscr]{eucal}
\theoremstyle{plain}
\usepackage{color}

\newtheorem{theorem}{Theorem}[section]
\newtheorem{lemma}[theorem]{Lemma}
\newtheorem{definition-theorem}[theorem]{Definition-Theorem}
\newtheorem{proposition}[theorem]{Proposition}
\newtheorem{corollary}[theorem]{Corollary}
\newtheorem{definition}[theorem]{Definition}
\newtheorem{example}[theorem]{Example}
\newtheorem{remark}[theorem]{Remark}
\newtheorem{conjecture}[theorem]{Conjecture}
\newtheorem{notation}[theorem]{Notation}

\newtheorem*{maintheorem*}{Main Theorem}
\newcommand \bth[1] { \begin{theorem}\label{t#1} }
\newcommand \ble[1] { \begin{lemma}\label{l#1} }

\newcommand \bpr[1] { \begin{proposition}\label{p#1} }
\newcommand \bco[1] { \begin{corollary}\label{c#1} }
\newcommand \bde[1] { \begin{definition}\label{d#1}\rm }
\newcommand \bex[1] { \begin{example}\label{e#1}\rm }
\newcommand \bre[1] { \begin{remark}\label{r#1}\rm }
\newcommand \bcj[1] { \begin{conjecture}\label{j#1}\rm }

\newcommand \bnota[1] { \begin{notation}\label{n#1}\rm }
\renewcommand {\eth} { \end{theorem} }
\newcommand {\ele} { \end{lemma} }

\newcommand {\epr} { \end{proposition} }
\newcommand {\eco} { \end{corollary} }
\newcommand {\ede} { \end{definition} }
\newcommand {\eex} { \end{example} }
\newcommand {\ere} { \end{remark} }
\newcommand {\ecj} { \end{conjecture} }

\newcommand {\enota} { \end{notation} }
\newcommand \thref[1]{Theorem \ref{t#1}}
\newcommand \leref[1]{Lemma \ref{l#1}}
\newcommand \prref[1]{Proposition \ref{p#1}}
\newcommand \coref[1]{Corollary \ref{c#1}}
\newcommand \cjref[1]{Conjecture \ref{j#1}}
\newcommand \deref[1]{Definition \ref{d#1}}
\newcommand \exref[1]{Example \ref{e#1}}
\newcommand \reref[1]{Remark \ref{r#1}}

\makeatletter
\newsavebox{\@brx}
\newcommand{\llangle}[1][]{\savebox{\@brx}{\(\m@th{#1\langle}\)}%
  \mathopen{\copy\@brx\kern-0.5\wd\@brx\usebox{\@brx}}}
\newcommand{\rrangle}[1][]{\savebox{\@brx}{\(\m@th{#1\rangle}\)}%
  \mathclose{\copy\@brx\kern-0.5\wd\@brx\usebox{\@brx}}}
\makeatother
\usepackage{hyperref}

\usepackage{tikz}
\usetikzlibrary{matrix}
\usepackage{graphicx,booktabs}
\usetikzlibrary{arrows.meta}
\usepackage{tikz-cd}

\DeclareMathOperator{\Ind}{Ind}

\DeclareMathOperator{\Tor}{Tor}
\DeclareMathOperator{\fp}{fp}

\DeclareMathOperator{\Spec}{Spec} 
\DeclareMathOperator{\Spech}{Spec^h}

 \DeclareMathOperator{\Proj}{Proj}

\DeclareMathOperator{\End}{End} 
 \DeclareMathOperator{\Hom}{Hom}

\DeclareMathOperator{\Mod}{{\sf Mod}}
\DeclareMathOperator{\modd}{{\sf mod}}

\DeclareMathOperator{\Ab}{{\sf Ab}}

\newcommand{\Moddash}[1]{\Mod \hspace{-0.3mm}\text{-} #1}
\newcommand{\moddash}[1]{\modd\hspace{-0.3mm}\text{-} #1}
\DeclareMathOperator{\Spc}{Spc}

\DeclareMathOperator{\colim}{colim}

\DeclareMathOperator{\Id}{Id}
\DeclareMathOperator{\ev}{ev}
\DeclareMathOperator{\coev}{coev}

\DeclareMathOperator{\cone}{cone}
\DeclareMathOperator{\h}{h}
\DeclareMathOperator{\Nil}{Nil}

\DeclareMathOperator{\cl}{cl}

\newcommand{\Spch}{\Spc^{\h}}

\newcommand{\mc}{\mathcal}

\newcommand{\id}{\operatorname{id}}

\newcommand{\bA}{\mathbf A}
\newcommand{\bT}{\mathbf T}

\newcommand{\bS}{\mathbf S}
\newcommand{\bK}{\mathbf K}
\newcommand{\bM}{\mathbf M}
\newcommand{\bC}{\mathbf C}
\newcommand{\bL}{\mathbf L}
\newcommand{\bP}{\mathbf P}
\newcommand{\bQ}{\mathbf Q}

\newcommand{\bI}{\mathbf I}
\newcommand{\bJ}{\mathbf J}

\newcommand{\Loc}{\operatorname{Loc}}
\newcommand{\XX}{\mathcal X}

\newcommand{\unit}{\ensuremath{\mathbf 1}}

\newcounter{listequation}

\numberwithin{equation}{section}

\begin{document}
\title[Homological Spectra of Monoidal Triangulated Categories]{The Homological Spectrum for \\
Monoidal Triangulated Categories}
\author[Daniel K. Nakano]{Daniel K. Nakano}
\address{Department of Mathematics \\
University of Georgia \\
Athens, GA 30602\\
U.S.A.}
\email{nakano@math.uga.edu}
%\thanks{}
\author[Kent B. Vashaw]{Kent B. Vashaw}
\address{
Department of Mathematics and Statistics \\
University at Albany, SUNY \\
Albany \\
NY 12222
\\
U.S.A.
}
\email{kvashaw@albany.edu}
%\thanks{}
\author[Milen T. Yakimov]{Milen T. Yakimov}
\address{
Department of Mathematics \\
Northeastern University \\
Boston, MA 02115 \\
U.S.A.
and International Center for Mathematical Sciences, Institute of Mathematics and Informatics \\
Bulgarian Academy of Sciences \\ 
Acad. G. Bonchev Str., Bl. 8 \\
Sofia 1113, Bulgaria
}
\email{m.yakimov@northeastern.edu}
\begin{abstract} The authors develop a notion of homological prime spectrum for an arbitrary monoidal triangulated category, ${\mathbf C}$. Unlike the symmetric case due to Balmer, the homological primes of ${\mathbf C}$ are not defined as the maximal Serre ideals of the small module category ${\sf mod}\text{-}{{\mathbf C}}$, or via a noncommutative ring theory inspired version of this construction. Instead, the authors work with an extended comparison map from the Serre prime spectrum $\operatorname{Spc} ({\sf mod}\text{-}{{\mathbf C}})$ to the Balmer spectrum $\operatorname{Spc}{\mathbf C}$, and select the maximal elements of each fiber to define the homological spectrum $\operatorname{Spc}^{\text{h}}{\mathbf C}$.

A surjective continuous homological comparison map $\operatorname{Spc}^{\text{h}}{\mathbf C} \to \operatorname{Spc}{\mathbf C}$ is constructed and used to formulate an extended Nerves of Steel Condition, stating that this map is a homeomorphism. It is believed that this condition holds for many naturally arising monoidal triangulated categories.

Nerves of Steel is shown to hold under stratification and uniformity conditions on ${\mathbf C}$. The proof is based on a general theorem giving an explicit description of the Balmer--Favi support of the pure-injectives associated to all Serre primes of ${\sf mod}\text{-}{{\mathbf C}}$. It is shown that the validity of Nerves of Steel carries over from ${\mathbf C}$ to a semidirect product ${\mathbf C} \rtimes G$ with an arbitrary group $G$, and when $G$ is infinite, this provides examples of monoidal triangulated categories satisfying Nerves of Steel and whose thick ideals are not centrally generated. Important cases in which the condition is verified include the stable module categories of the coordinate rings of all finite group schemes and the Benson--Witherspoon Hopf algebras.  
\end{abstract}
\maketitle
\section{Introduction}
\label{Introduction}

The Balmer spectrum \cite{Balmer0} of a symmetric monoidal triangulated category $\bC$ is a topological space $\Spc \bC$ (equipped with a locally ringed structure) that captures many fundamental properties of 
$\bC$. For example, one of Balmer's results \cite{Balmer0} gives an explicit classification of all thick ideals of $\bC$ in terms of $\Spc \bC$. However, the description of $\Spc \bC$ is very difficult and explicit results are only known in several cases. In all of such cases one resorts on an explicit support map 
\[
\sigma : \bC \to \XX_{\cl}(Y)
\]
with the tensor product property
\[
\sigma (A \otimes B) = \sigma (A) \cap \sigma(B), \quad \forall A, B \in \bC,
\]
see \deref{sup}. Here $Y$ is a topological space, $\XX_{\cl}(Y)$ is the collection of closed subsets of $Y$. Often times, the collection of objects of a category will be denoted by the same symbol. 

In order to address the tensor product property for finite groups, one relies on Carlson's support \cite{Carlson} via rank varieties. For finite group schemes, a crucial role is played by the Friedlander--Pevtsova $\pi$ point support \cite{FriedlanderPevtsova}. The classification results 
of Hopkins \cite{Hopkins}, Neeman \cite{Neeman}, and Thomason \cite{Thomason} for the Balmer spectrum of the perfect derived category of a scheme is based on the classical sheaf support. Unlike these important support maps that work in concrete settings, there is no known argument that proves that the generally defined cohomological support map has the tensor product property. 

For all symmetric monoidal triangulated categories $\bC$,
Balmer constructed a support map with the tensor product property, the {\em{homological support}}, that conjecturally can be used in place of the specific support maps in the above cases \cite{Balmer-support-big,Balmer2020c}. The construction is based on the Freyd envelope of the small module category $\moddash{\bC}$. The topological space consists of maximal Serre ideals of $\moddash{\bC}$; it is termed the {\em{homological spectrum}} of $\bC$ and denoted by $\Spch \bC$. There is a comparison map that is surjective under natural rigidity conditions  
\[
\Spch \bC \to \Spc \bC
\]
that intertwines the homological support map and the universal support map based on the Balmer spectrum of $\bC$.

The notion of tensor-triangular fields was introduced and studied by Balmer, Krause, and Stevenson \cite{BKS2019}. Balmer proved an important nilpotence theorem for maps in $\bC$ using elements in $\Spch \bC$ and 
homological residue fields \cite[1.1 Theorem]{Balmer2020c}. Furthermore, he posed what has become known as the Nerves of Steel Conjecture that the above comparison map is a homeomorphism \cite{Balmer-support-big}. Recently, counterexamples to the Nerves of Steel Conjecture have been constructed by Barthel--Hyslop--Ramzi \cite{BHR}, and, separately, by Coulembier \cite{Coulembier}. It is still believed that Nerves of Steel holds for many naturally-arising tensor-triangulated categories.

Balmer and Cameron \cite{BalmerCameron2021} explicitly described homological residue fields for tensor triangulated categories which are of wide interest. Recently, the homological spectrum of certain classes of Lie superalgebra representations was investigated in \cite{Hamil-Nakano,Hamil} where Nerves of Steel was proved in super type A.
 Nerves of Steel is closely related to a certain exact nilpotence condition for local tensor-triangulated categories. This exact nilpotence condition has recently been verified by Hyslop for a wide class of $\mathbb{E}_{\infty}$-rings \cite{Hyslop,Hyslop2026}.

For arbitrary monoidal triangulated category $\bC$ that need not be symmetric, the authors have developed the foundations of 
{\em Noncommutative Tensor Triangular Geometry}, which is the study of the space $\Spc \bC$ in \cite{NVY1,NVY2,NVY3,NVY4}. Many important representation theoretic settings involve monoidal triangulated categories that are not even braided, such as the stable categories of non-quasitraingular Hopf algebras, stable categories of modules of vertex operator algebras, and finite tensor categories \cite{EGNO}. These categories play a prominent role in representation theory, mathematical physics and topological quantum computing, and are the subject of intense current research. 

The idea of thick prime ideals of a monoidal triangulated category $\bC$ was inspired from noncommutative ring theory by translating the defining property of the latter to a triangulated setting \cite{BKS,NVY1}. Even though this might have seemed naive at first sight, it turned out to 
highly effective for handling a wider variety of categories. For instance, under mild conditions $\Spc \bC$ naturally classifies all thick ideals of $\bC$, as shown in \cite{NVY4,Rowe}. Moreover, the differences between
the noncommutative Balmer spectrum allowed the authors 
to detect the tensor product property of support maps and to formulate a precise conjecture on the description of the Balmer spectrum for an arbitrary finite tensor categories. See  \cite{NVY3} and \cite{NVY4}. 

Given the recent counterexamples to Nerves of Steel, developing the notion of homological primes with the associate topological space $\Spch \bC$ for monoidal triangulated categories is an important and challenging problem in the field. 
This is at the heart of what we propose in this paper. For the setting, we work with rigidly-compactly generated monoidal triangulated categories $\bK$ with compact part denoted by $\bK^c$ (see \deref{main-def}), and consider the small module category $\moddash{\bK^c}$ with the Day convolution product, which is generally only right exact (see Section \ref{3} for details). In our construction, we focus on all Serre ideals of $\moddash{\bK^c}$ and do not pick up just the maximal ones since those could be too few. Instead, we look at all Serre prime ideals. The corresponding space $\Spc (\moddash{\bK^c})$ was considered in \cite{Vashaw-Yakimov} in relation to the categorification of Richardson varieties. In this setting, 
we prove that the inverse images under the Yoneda embedding 
$\h : \bK^c \to \moddash{\bK^c}$ define a continuous map 
\[
\widetilde{\phi} : \Spc(\moddash{\bK^c}) \to \Spc \bK^c, 
\quad \widetilde{\phi}(\bS) := \h^{-1} (\bS). 
\]

The definition presented below will be a key idea that is used throughout the paper.
\medskip

\noindent
{\bf{Definition.}} 
\hfill
\begin{enumerate}
\item[(a)] The {\em{homological primes}} of a
monoidal triangulated category $\bK^c$ are the maximal (with respect to inclusion) elements in each fiber of 
$\widetilde{\phi}$.
\item[(b)] The {\em{homological spectrum}} $\Spch \bK^c$ of $\bK^c$ is the space of homological primes with the topology where closed sets are defined to be arbitrary intersections of the base of closed sets
\begin{equation}
    \label{Vhdef}
V^{\h}(A):=\{\bS \in \Spch \bK^c : \h(A) \not \in \bS\}
\end{equation}
for $A \in \bK^c$.
\item[(c)] The {\em{homological support map}} $V^{\h} : \bK^c \to \XX_{\cl}(\Spch \bK^c)$ is the map given by \eqref{Vhdef}.
\end{enumerate}
\medskip

The formula defining $\widetilde{\phi}$ produces a map from the set of thick ideals of $\bK^c$ to the set of thick ideals of $\moddash{\bK^c}$. A homological prime of $\bK^c$ 
can be equivalently defined as a maximal element in a fiber of this extension of $\widetilde{\phi}$ over a prime ideal of $\bK^c$. Note that this second definition does not refer to Serre primes of $\moddash{\bK^c}$ at all. 

In order to connect the aforementioned definitions of homological primes to maximal Serre ideals, one needs the following condition: 
\smallskip

\noindent 
{\em{$\bK$ satisfies the {\sf Central Generation Hypothesis} if all thick ideals of $\bK^c$ are generated by collections of central elements, which are elements $A \in \bK^c$ such that $A \otimes B \cong B \otimes A$ for every object $B \in \bK^c$ in a possibly nonfunctorial way.}}
\medskip

\noindent 
{\bf{Proposition.}} {\em{If $\bK$ satisfies the Central Generation Hypothesis, then the homological primes of $\bK^c$ are precisely the maximal Serre ideals of $\moddash{\bK^c}$.}}
\medskip

In particular, if $\bK^c$ is braided, our definition recovers Balmer's notion of homological primes \cite{Balmer-support-big,Balmer2020c}. Furthermore, we demonstrate that in all cases a homological prime $\bS$ has the weaker maximality property that for every central element $A \in \bK^c$ such that $\h(A) \notin \bS$, 
the minimal Serre ideal generated by $\bS$ and $\h(A)$ is 
all of $\moddash{\bK^c}$. 
\smallskip

The important property of the homological support map is that it has the noncommutative tensor product property as defined by the authors 
in \cite{NVY1}:
    \[
    \bigcup_{C \in \bK^c} V^{\h} (A   \otimes C \otimes B) = V^{\h} (A) \cap V^{\h} (B), \quad \forall A, B \in \bK^c.
    \]
The restriction of $\widetilde{\phi}$ to $\Spch \bK^c$ will 
be called {\em{homological comparison map}} and will be denoted by 
\[
\phi : \Spch \bK^c \to \Spc \bK^c.
\]
It intertwines the homological support map and the Balmer support map $V^B: \bK^c \to \Spc \bK^c$ given by $V^B(A) = \{ \bP \in \Spc \bK^c : A \notin \bP \}$ in the sense that 
$\phi^{-1} (V^B(A))= V^{\h}(A)$, for all $A \in \bK^c$. 
Such a comparison map exists for all support maps on $\bK^c$ 
that satisfy the tensor product property \cite{BKS,NVY1}. Due to the surjectivity of $\widetilde{\phi}$, 
the homological comparison map is a surjective continuous map. 
With our setup, one can formulate a general version of Nerves of 
Steel. In view of the counterexamples of Barthel--Hyslop--Ramzi and Coulembier, we refer to the Nerves of Steel Condition rather than Conjecture. 
\medskip

\noindent
{\bf{Condition.}} [{\em{Nerves of Steel\ }}]
{\em{For a rigidly-compactly generated monoidal triangulated category $\bK$, the homological comparison map 
\[
\phi : \Spch \bK^c \to \Spc \bK^c
\]
is a homeomorphism.}}
\medskip

The above framework for homological primes and homological comparison map is developed in Sections \ref{sec:hom-prim} and \ref{S:supportmap}. 

The second goal of this paper to develop methods towards verifying Nerves of Steel for a given monoidal triangulated category. Krause \cite{Krause2000} proved that the restricted Yoneda functor gives a bijection between pure-injective objects of $\bK$ and injective objects of $\Moddash{\bK^c}$. 

Firstly, in Section \ref{6} we provide a direct extension of the results of \cite{BKS2019}, to show that for every Serre ideal $\bS$, there exists a canonical pure-injective object $E_\bS$ of $\bK$ such that 
\[
\ker( \h(E_\bS) \otimes -) = \ker(- \otimes \h(E_\bS))
\]
equals the minimal localizing subcategory of $\Moddash{\bK^c}$ containing $\bS$ (which is necessarily an ideal of $\Moddash{\bK^c}$). In order to apply this property towards the conjecture, one needs to know the Balmer--Favi \cite{Balmer-Favi} support of the pure-injectives associated to 
all Serre primes of $\bK^c$. Note that the Balmer--Favi support was extended to the noncommutative situation in \cite{CaiVashaw}, see Section \ref{7} for details.
\medskip

\noindent 
{\bf{Definition.}}
We say that $\bQ \in \Spc(\bK^c)$ {\em{lives over}} a Serre prime ideal $\bS$ of $\moddash{\bK^c}$ if there exists a Serre prime $\bS'$ of $\moddash{\bK^c}$ containing $\bS$ such that $\phi(\bS')=\bQ$. 
\medskip

For example, note that, $\phi(\bS) \in \Spc (\bK^c)$ lives over $\bS$. 
Our first main result is as follows:
\medskip

\noindent
{\bf{Theorem A.}} {\em{Assume that $\bK$ is a rigidly-compactly generated monoidal category whose Balmer spectrum $\Spc \bK^c$ is Noetherian. Let $\bS$ be a Serre prime ideal of 
$\moddash{\bK^c}$ and $E_{\bS}$ be the corresponding pure-injective object of $\bK$. Then the Balmer--Favi support of  $E_{\bS}$ is equal to}}
\[
V^{BF} (E_{\bS})=\{\phi(\bS')\in \Spc \bK^c : \bS' \in \Spch(\bK^c) \; \; \text{such that} \; \bS' \supseteq \bS \} \subseteq \Spc \bK^c.
\]
In other words, the Balmer--Favi support of $E_{\bS}$ consists of all prime ideals of $\bK^c$ living over $\bS$. 
\smallskip

The theorem is of independent interest since the Balmer--Favi support has many additional applications. The theorem is new even in the symmetric case and is a far-reaching extension of \cite[Lemma 3.7]{BHS2023} treating the case
of maximal Serre primes of $\moddash{\bK^c}$, in which case 
the support is a singleton $V^{BF}(E_\bS) = \{\phi(\bS)\}$. 
\medskip

The following hypothesis will be important in understanding  Nerves of Steel: 
\smallskip

\noindent 
{\em{$\bK$ satisfies the {\sf Uniformity Hypothesis} if for all $\bS_1, \bS_2 \in \Spch \bK^c$ in the same fiber of $\phi$, i.e. $\phi(\bS_1) = \phi(\bS_2)$, the collections of primes of $\bK^c$ living over $\bS_1$ and $\bS_2$, respectively, coincide.}}
\smallskip

For example, every category $\bK$ satisfying the Central Generation Hypothesis, satisfies the Uniformity Hypothesis. In particular, all braided categories $\bK$ satisfy the Uniformity Hypothesis. Our second main result gives conditions that imply the Nerves of Steel Condition. 
\medskip

\noindent
{\bf{Theorem B.}} {\em{Assume that $\bK$ is a rigidly-compactly generated monoidal triangulated category, which is stratified by $\Spc \bK^c$, satisfies the Uniformity Hypothesis, and that the Balmer spectrum $\Spc \bK^c$ is Noetherian. Then $\bK$ satisfies Nerves of Steel.}}
\medskip

Theorems A and B are proved in Section \ref{7}. On the basis of the second theorem we can analyze the semidirect products of monoidal categories and (generally infinite) groups $G$ as follows:
\medskip

\noindent
{\bf{Theorem C.}} {\em{Assume that $\bK$ is a rigidly-compactly generated monoidal triangulated category satisfying Nerves of Steel, the Central Generation Hypothesis, and that the Balmer spectrum $\Spc \bK^c$ is Noetherian. Let $G$ be a group acting on $\bK$ by monoidal triangulated autoequivalences. Then $\bK \rtimes G$ satisfies Nerves of Steel.}}
\medskip

For infinite groups $G$,  Theorem C presents examples of categories of the form $\bK \rtimes G$ which satisfy Nerves of Steel and do not satisfy the Central Generation Hypothesis. This provides examples of categories satisfying the conjecture that are as far from braided as one could imagine.
On the other hand, one can show that if $\bK$ is symmetric and $G$ is finite, then the category $\bK \rtimes G$ satisfies the Central Generation Hypothesis. 
Theorem C is proved in Section \ref{8}.

Theorem C can be used to prove Nerves of Steel for important monoidal triangulated categories. This is done in Section \ref{9}, where 
we prove the condition for the stable module categories of the coordinate rings of all finite groups schemes and Benson--Witherspoon Hopf algebras \cite{Benson-Witherspoon}. In Section \ref{9} we also put in perspective Nerves of Steel for the stable categories of all finite tensor categories \cite{Etingof-Ostrik} with the previous works on the Balmer spectra of those categories; we conjecture that for these specific categories, Nerves of Steel always holds. 
This conjecture fits in a broader conjecture that both the homological and cohomological comparison maps for the Balmer spectra of the stable categories of all finite tensor categories are homeomorphisms. The homological spectrum and its canonical support map are the natural candidates for a support map with the noncommutative tensor product property \cite{NVY1} which is to replace the rank rank support of Carlson \cite{Carlson} and the $\pi$-support of Friedlander--Pevtsova \cite{FriedlanderPevtsova} for the much wider generality of finite tensor categories. 
\medskip

\noindent
{\bf{Acknowledgments.}} 
The authors thank Kevin Coulembier for alerting them to the counterexample in \cite{Coulembier}.

The research of D.K.N was supported in part by NSF grant DMS-2401184. The research of  K.B.V. was supported in part by NSF Postdoctoral Fellowship DMS-2103272 and by an AMS-Simons Travel Grant. The research of M.T.Y. was supported in part by NSF grant DMS–2601486, the Simons Foundation, grant SFI-MPS-T-Institutes-00007697, and the Bulgarian Ministry of Education and Science grant DO1-239/10.12.2024 and Science Fund grant KP-06-N92/5.

\section{The noncommutative Balmer spectrum}
\noindent
\bde{main-def}
A {\em{rigidly-compactly generated monoidal triangulated category}} is a triangulated category $\bK$ with compact part $\bK^c$ (the full subcategory of objects $C \in \bK$ such that $\Hom_{\bK}(C, -)$ commutes with set arbitrary coproducts) such that the following hold:
\begin{enumerate}
    \item[(i)] $\bK$ admits arbitrary set-indexed coproducts;  
    \item[(ii)] $\bK^c$ generates $\bK$ (for every object $M$ of $\bK$ there exists an object $M$ of $\bK^c$ such that $\Hom_{\bK}(M,C) \neq 0$);
    \item[(iii)] $\bK$ is a monoidal category with unit object $\unit \in \bK^c$, such that the monoidal product is biexact and commutes with arbitrary coproducts;  
    \item[(iv)] $\bK^c$ is a rigid monoidal subcategory of $\bK$.
\end{enumerate}
\ede

Recall that a monoidal category is called rigid if every object $A$ has left and right duals $A^*$ and ${}^* A$ with appropriate evaluation and coevaluation maps, \cite[Definitions 2.10.1 and 2.10.2]{EGNO}. The opposite implication to the one in condition (4) always holds: If $M \in \bK$ is a left or right dualizable then it is compact. Indeed, in the case when $M$ has a left dual $M^*$ and by 
\cite[Proposition 2.10.8(a)]{EGNO}, we have
\[
\Hom_\bK(M, -) \cong \Hom_{\bK} (\unit \otimes M, - ) 
\cong \Hom_{\bK} ( \unit, M^* \otimes -). 
\]
Since $\unit \in \bK^c$ and the functor $M^*\otimes -$ commutes with arbitrary coproducts, the functor $\Hom_\bK(M, -)$ has the same property. 

{\em{Symmetric rigidly-compactly generated monoidal triangulated categories}} have been much studied, see e.g.,\cite{Balmer-support-big,BKS2019,BKS2020}. They are also referred to as {\em{big tt categories}}.  
\medskip

\noindent
{\bf{Assumption.}} {\em{Throughout the paper we will work with rigidly-compactly generated monoidal triangulated categories $\bK$.}}
\medskip

On the side of the small category $\bK^c$, we will work with the following constructions:

\bde{def-ideals-small}
\hfill
\begin{itemize}
\item[(a)] A {\em{thick ideal}} of $\bK^c$ is a full triangulated subcategory $\bI$ of $\bK^c$ that is closed under direct summands such that both objects $A \otimes B$ and $B \otimes A$ are in $\bI$ for all $A \in \bI$ and $B \in \bK^c$.
For a set of objects $\bS$ of $\bK^c$, we will denote by $\langle \bS \rangle$ the (unique) minimal thick ideal of $\bK^c$ containing $\bS$. 
\item[(b)] A proper thick ideal $\bP$ of $\bK^c$ is called {\em{prime}} 
if $\bI \otimes \bJ \subseteq \bP$ implies that one of $\bI$ or $\bJ$ is contained in $\bP$, for all thick ideals $\bI$ and $\bJ$ of $\bK^c$. 
This is equivalent to the property that for all $A, B \in \bK^c$, 
\[
A \otimes C \otimes B \in \bP, \; \forall C \in \bK^c \quad \Rightarrow A \in \bP \; \mbox{or} \; B \in \bP,
\]
see \cite[Theorem 1.2.1(a)]{NVY1}.
\item[(c)] A proper thick ideal $\bP$ of $\bK^c$ is called {\em{completely prime}} 
if $A \otimes B\in \bP$ implies that one of $A$ or $B$ is an object in $\bP$. 
\end{itemize} 
\ede
{\em{The noncommutative Balmer spectrum}} of $\bK^c$, denoted by 
$\Spc \bK^c$, is the collection of all prime ideals of $\bK^c$, equipped with the Zariski topology where closed sets are defined 
to be arbitrary intersections of the base of closed sets
\begin{equation}
\label{Balmer-sup}
V^B (A) = \{ \bP \in \Spc \bK^c : A \not \in \bP \}.
\end{equation}
This space was introduced by Balmer \cite{Balmer0} in the symmetric case and by Buan--Krause--Solberg \cite{BKS} in general. 
The authors extensively developed the foundations of this spectrum and made important calculations in \cite{NVY1,NVY2,NVY3}. Although the notion of completely prime ideal appears more attractive at first sight because it directly mimics the classical situation developed in \cite{Balmer0}, it is more restrictive, and equality of the sets of prime and completely prime ideals only occurs when the universal support map $A \mapsto V(A)$ has the tensor product property $V(A \otimes B) = V(A) \otimes V(B)$ for all objects $A$ and $B$ of $\bK^c$, \cite[Theorem 3.1.1]{NVY2}.

Under mild conditions, the set of all thick ideals of $\bK^c$ can be reconstructed from the noncommutative Balmer spectrum 
$\Spc \bK^c$, which is one of the important roles of the latter.
The following theorem was proved in the symmetric case by Balmer and in general in \cite[Proposition A.7.1]{NVY4}. Under the additional assumption that $\Spc \bK^c$ is Noetherian, it was also obtained by Rowe 
\cite{Rowe}. 

\bth{thm-bij} If $\bK^c$ is generated by a single object as a thick subcategory then there is a bijection between the collection of thick ideals of $\bK^c$ and the collection of Thomason subsets of $\Spc \bK^c$ (unions of closed subsets whose complements are quasicompact). 
\eth

For a topological space $Y$, denote by $\XX_{\cl}(Y)$ the collection of closed subsets of $X$. 
\bde{sup}
\hfill 
\begin{enumerate}
    \item[(a)] A support map for $\bK^c$ is a map $\sigma : \bK^c \to \XX_{\cl}(Y)$ for a topological space $X$ satisfying the following conditions 
    \begin{enumerate}
\item[(i)] $\sigma(0)=\varnothing $ and $ \sigma(1)= Y$;
\item[(ii)] $\sigma(A\oplus B)=\sigma(A)\cup \sigma(B)$, $\forall A, B \in \bK^c$; 
\item[(iii)] $\sigma(\Sigma A)=\sigma(A)$, $\forall A \in \bK^c$; 
\item[(iv)] If $A \to B \to C \to \Sigma A$ is a distinguished triangle, then $\sigma(A) \subseteq \sigma(B) \cup \sigma(C)$:
\end{enumerate}
    \item[(b)] We say that a support map $\sigma : \bK^c \to \XX_{\cl}(Y)$ has the noncommutative tensor product property if 
    \[
    \bigcup_{C \in \bK^c} \sigma (A   \otimes C \otimes B) = \sigma (A) \cap \sigma(B), \quad \forall A, B \in \bK^c
    \]
\end{enumerate}
\ede
The Balmer support map 
\begin{equation}
\label{eq:Balmer-sup}
V^B : \bK^c \to \Spc (\bK^c)
\end{equation}
is given by \eqref{Balmer-sup}. It is a support map in the sense of \deref{sup} and has the noncommutative tensor product property. The following universal property of the Balmer support was proved in \cite{Balmer0} in the symmetric case and \cite[Theorem 4.5.1]{NVY1} in general. 

\bpr{univ} The Balmer support is a final object in the category of all support maps having the noncommutative tensor product property: for every such support map $\sigma : \bK^c \to \XX_{\cl}(Y)$  there is a unique continuous map $f_\sigma: Y \to \Spc \bK^c$ such that 
\[
\sigma(A)= f_\sigma^{-1}(V^B(A))
\]
for all $A \in \bK^c$. 
\epr
\bre{restr} If $\sigma : \bK^c \to \XX_{\cl}(Y)$ is a support map and $Y'$ is a subset of $Y$ equipped with the induced topology, then 
\[
\sigma' : \bK^c \to Y', 
\quad 
\sigma'(A):= \sigma(A) \cap Y', \; \; \forall A \in \bK^c
\]
is also a support map. If in addition, $\sigma$
has the tensor product property, then $\sigma'$ has it too and the associated map $f_\sigma: Y' \to \Spc \bK^c$
from \prref{univ} is the restriction
\[
f_{\sigma'} = f_{\sigma}|_{Y'}. 
\]
\ere

Finally, on the side of the big category $\bK$, we will work with the following constructions:

\bde{def-loc-big}
\hfill
\begin{enumerate}
    \item[(a)] A {\em{localizing subcategory}} of $\bK$ is a triangulated category that is closed under arbitrary coproducts. 
    \item[(b)] A {\em{localizing ideal}} $\bI$ of $\bK$ is a localizing subcategory such that $A \otimes B, B \otimes A \in \bI$ for all $A \in \bI$ and $B \in \bK$. For a collection of object $\bS$ of $\bK$ we will denote by $\Loc(\bS)$ the (unique) minimal localizing ideal of $\bK$ containing $\bS$. 
\end{enumerate}
\ede
\section{The module category of a monoidal triangulated category}
\label{3}
Let $\bK$ be a rigidly-compactly generated triangulated category with compact part $\bK^c$. Denote by $\Moddash{\bK^c}$ the {\em{module}} 
(or {\em{functor}}) {\em{category}} of $\bK^c$, which is 
the (abelian) category of contravariant additive functors $\bK^c \to \Ab$, where $\Ab$ denotes the category of abelian groups. We have the {\em{Yoneda embedding}}
\[
\h : \bK^c \hookrightarrow \Moddash{\bK^c}
\]
and its extension, 
\[
\h: \bK \to \Moddash{\bK^c}, \quad 
\quad \h(A):= \Hom_{\bK} (-,A),
\]
called the {\em{restricted Yoneda functor}}. 

The category $\Moddash{\bK^c}$ is a {\em{Grothendieck category}}, that is, it is an abelian category admitting arbitrary set-indexed coproducts (hence arbitrary set-indexed colimits), direct limits of exact sequences are exact (direct limits are always right exact, so the salient point is left exactness), and it has a generator. Here a generator means an object $G$ such that every object is a quotient of some (possibly infinite) direct sum of $G$'s. 

The {\em{Freyd envelope}} of the {\em{small module category}} of $\bK^c$ is defined to be
the full subcategory of {\em{finitely presented objects}}
\begin{equation}
\label{eq:modMod}
\moddash{\bK^c}:= \left(\Moddash{\bK^c}\right)^{\fp}. 
\end{equation}
Recall that an object $M \in \Moddash{\bK^c}$ is finitely presented if 
\[
M \mapsto \Hom_{\Moddash{\bK^c}}(M, -) 
\]
commutes with filtered colimits. 
The Yoneda embedding lands in the Freyd envelope: 
$\bK^c \hookrightarrow \moddash{\bK^c}$. Every object of $\Moddash{\bK^c}$ is a filtered colimit of objects of 
$\moddash{\bK^c}$.

The suspension $\Sigma : \bK^c \to \bK^c$ induces a suspension functor $\Sigma : \Moddash{\bK^c} \to \Moddash{\bK^c}$ given by on objects by 
\begin{equation}
\label{Susp}
\Sigma (F)(A) := F(\Sigma^{-1}(A)) 
\end{equation}
for all $A \in \bK^c$, $F \in \Moddash{\bK^c}$. 

For our setting, it is important to describe various subcategories and ideals in $\moddash{\bK^c}$ and $\Moddash{\bK^c}$. 

\bde{subcatsMod}
\hfill
\begin{enumerate}
\item[(a)] A Serre ideal of $\moddash{\bK^c}$ is a 
Serre subcategory $\bI$ of $\moddash{\bK^c}$ that is closed under left and right tensoring with arbitrary objects of $\moddash{\bK^c}$. 
 \item[(b)] For a collection of objects $\bS$ of $\moddash{\bK^c}$, denote by $\langle \bS \rangle$
the (unique) Serre ideal of $\moddash{\bK^c}$ containing $\bS$.
\item[(c)] A localizing subcategory of $\Moddash{\bK^c}$ is a Serre subcategory $\bS$ of $\Moddash{\bK^c}$ that is closed under arbitrary coproducts (or equivalently under arbitrary colimits). A localizing ideal of $\Moddash{\bK^c}$ is a localizing subcategory that is closed under left and right tensoring with arbitrary objects of $\Moddash{\bK^c}$. 
\end{enumerate}
\ede

Some additional notions will be introduced in Section \ref{sec:hom-prim}. For brevity we will denote 
\[
\hat{A} := \h(A), \quad \forall A \in \bK.
\]
Every Serre ideal of $\moddash{\bK^c}$ is invariant under the suspension functor \eqref{Susp} since 
\[
\Sigma^k(C) \cong \Sigma^n(\hat{\unit}) \otimes C, 
\quad 
\forall C \in \moddash{\bK^c}. 
\]
The restricted Yoneda functor $\bK \stackrel{\h}{\rightarrow} \Moddash{\bK^c}$ is a homological functor (sending distinguished triangles to long exact sequences) and preserves arbitrary coproducts. In addition, it is a universal functor with these properties, which means that every
homological functor that preserves arbitrary coproducts $\bK \stackrel{f}{\rightarrow} \bC$ coincides with the composition of $h$ and an exact functor that preserves coproducts $\Moddash{\bK^c} \stackrel{\hat{f}}{\rightarrow} \bC$. There is an explicit formula for the functor $\hat{f}$, \cite[Eq. (A.11)]{BKS2020}, which we review next. 

Firstly, recall that for an object $M \in \Moddash{\bK^c}$, the slice category $[\bK^c \stackrel{\h}{\to}M]$ has objects given by pairs $(\hat{A}, \psi \in \Hom_{\Moddash{\bK^c}}(\hat{A}, M))$ and morphisms $\hat{\theta} : (\hat{A}_1, \psi_1) \to (\hat{A}_1, \psi_2)$ given by $\theta \in \Hom_{\bK} (A_1, A_2)$ such that $\psi_2 = \psi_1 \theta$. Every object of $M \in \Moddash{\bK^c}$ is isomorphic to a colimit of objects in the image of the Yoneda embedding:
\begin{equation}
    \label{eq:colimM}
M = \underset{(A, \hat{A} \to M) \in [\bK^c \stackrel{\h}{\to}M]}{\colim} \hat{A}.
\end{equation}
Since the functor $\hat{f}$ is exact and commutes with arbitrary coproducts, it commutes with arbitrary colimits. 
Because $\hat{f}(A) = f(h(A))$, $\forall A \in \bK^c$, 
\[
\hat{f}(M) = \underset{(A, \hat{A} \to M) \in [\bK^c \stackrel{\h}{\to}M]}{\colim} f(A), \quad \forall M \in \Moddash{\bK^c}. 
\]

For $A \in \bK$, consider the homological, coproduct preserving functors from $\bK$ to $\Moddash{\bK^c}$:
\[
r_A(-) := \h(A \otimes -)
\quad \mbox{and} \quad 
l_A(-) := \h(- \otimes A). 
\]
Denote the corresponding exact, coproduct preserving endofunctors of $\Moddash{\bK^c}$ 
by $\hat{r}_A$ and $\hat{l}_A$, respectively. 

\ble{yoneda-tensor}
\hfill
\begin{enumerate}
\item[(a)] The module category $\Moddash{\bK^c}$ is monoidal with a right exact tensor functor, called the Day convolution product, defined on objects by 
\[
M \otimes N := 
\underset{(A, \hat{A} \to M) \in [\bK^c \stackrel{\h}{\to} M]}{\colim}
\; \; \; 
\underset{(B, \hat{B} \to N) \in [\bK^c \stackrel{\h}{\to} N]}{\colim}
\widehat{A \otimes B},  
\]
and unit object $\hat{\unit}$. The
small module category $\moddash{\bK^c}$ is a monoidal subcategory. 
\item[(b)] The restricted Yoneda functor $\h : \bK \to \Moddash{\bK^c}$ is a monoidal functor. 
\item[(c)] If $A \in \bK$, then the functors $\hat A \otimes -$ and $- \otimes \hat A$ coincide with $\hat{l}_A(-)$ and $\hat{r}_A(-)$, respectively. In particular, $\hat{A}$ is flat, i.e., $\hat A \otimes -$ and $- \otimes \hat A$ are exact endofunctors of $\Moddash{\bK^c}$. 
\item[(d)] The Yoneda embedding $\h : \bK^c \hookrightarrow \moddash{\bK^c}$ intertwines the Verdier and Gabriel--Serre quotients:
one has the monoidal abelian equivalence
\[
\moddash{\bK^c}/\langle \h(\bI)\rangle \cong \moddash{(\bK^c/\bI)}
\]
for all thick ideals $\bI$ of $\bK^c$. (Note that $\bK^c/\bI$ is a rigid monoidal triangulated category.) 
\end{enumerate}
\ele

The aforementioned lemma can be justified as follows. Part (a) of the lemma is well known. Proofs of parts (b) and (c) can be found in the appendix of
\cite{BKS2020}. The results in \cite{BKS2020} are stated for the case when $(\bK, \otimes)$ is a symmetric monoidal category, but the proofs do not use the assumption.  

A proof of part (d) is given in \cite[Proposition 2.9]{Balmer2020c}. Once again the setting of \cite[Proposition 2.9]{Balmer2020c} is phrased in the symmetric setting, but the same proof holds in general without any changes. 
\section{Homological primes}
\label{sec:hom-prim}
%%%%%%%%%%%%%%%%%%%%%%%%%%%%%%%%%%%%%%%%%%%%%%%%%

Throughout, let $\bK$ be a rigidly-compactly generated monoidal triangulated category with compact part $\bK^c$.

\bde{serre-prime}
\cite[Definition 6.1]{Vashaw-Yakimov}
A {\em{Serre prime}} in an abelian monoidal category $\bC$ is a Serre ideal $\bS$ of $\bC$ (recall \deref{subcatsMod}(b)) such that if $A \otimes \bC \otimes B \subseteq \bS$, for any $A$ or $B$ in $\bC$, then $A$ or $B$ is in $\bS$.
\ede
The case of biexect tensor product was considered in \cite{Vashaw-Yakimov}. Here, by the nature of $\moddash{\bK^c}$, we consider the more general case 
of a right exact tenor product in each argument. 

Denote the Serre monoidal spectrum $\Spc (\moddash{\bK^c})$ to the topological space consisting of Serre prime ideals of $\moddash{\bK^c}$ with the topology where closed sets are defined to be arbitrary intersections of the base of closed sets
\begin{equation}
\label{wtVh}
\widetilde{V}^{\h}(A):=\{\bS \in \Spc(\moddash{\bK^c}) : \hat A \not \in \bS\}
\end{equation}
for $A \in \bK^c$.

\ble{inv-im}
Let $\bS$ be a Serre prime of $\moddash{\bK^c}$. Then $\h^{-1}(\bS)$ is a prime ideal of $\bK$, i.e., an element of $\Spc \bK^c$.
\ele

\begin{proof}
Suppose $A \otimes \bK^c \otimes B \subseteq \h^{-1}(\bS)$, in other words, $\hat A \otimes \hat C \otimes \hat B \subseteq \bS$ for any object $C \in \bK^c$. Let $F \in \moddash{\bK^c}$. Pick $D \in \bK$ with $\hat D \twoheadrightarrow F$ in $\moddash{\bK^c}$. By \leref{yoneda-tensor}(c) we have 
\[
\hat A \otimes \hat D \otimes \hat B \twoheadrightarrow \hat A \otimes F \otimes \hat B.
\]
Since the left hand side is in $\bS$, and $\bS$ is a Serre subcategory, the right hand side is in $\bS$ as well. It follows that $\hat A \otimes \moddash{\bK^c} \otimes \hat B \subseteq \bS$. Since $\bS$ is prime, either $A$ or $B$ is in $\h^{-1}(\bS)$. 
\end{proof}

\leref{inv-im} produces the map
\begin{equation}
\label{eq:phi}
\widetilde{\phi} : \Spc(\moddash{\bK^c}) \to \Spc \bK^c, \quad \widetilde{\phi}(\bS) := \h^{-1}(\bS). 
\end{equation}

The foundational definition in this paper is the following.

\bde{homol-prime}
A {\em{homological prime}} for $\bK^c$ is a Serre prime ideal $\bS$ of $\moddash{\bK^c}$ which is maximal in the fiber of $\widetilde{\phi}(\bS)$ under the map \eqref{eq:phi}. 

In explicit terms, a Serre prime ideal $\bS$ of $\moddash{\bK^c}$ is a homological prime if it is a maximal element (with respect to inclusion) in the collection of Serre primes, $\widetilde{\phi}^{-1}( \widetilde{\phi}(\bS))$, of $\moddash{\bK^c}$. 
\ede

\bre{homol-prime-equiv}
    Equivalently, a homological prime may be defined as a Serre ideal $\bS$ such that $\bS$ is maximal in the fiber of $\widetilde{\phi}(\bS)$, and $\widetilde{\phi}(\bS) \in \Spc \bK^c$. This follows at once from the fact that the Yoneda embedding $\h : \bK^c \hookrightarrow \moddash{\bK^c}$
is a monoidal functor, \leref{yoneda-tensor}(b).
\ere

In the case that $\bK$ is braided, it is a well-known result that all homological primes are maximal in $\moddash{\bK^c}$. We have the following more involved version of this fact in the noncommutative situation. 

\ble{homol-prime-max}
Let $\bS$ be a homological prime in $\moddash{\bK^c}$. If $A$ is an object which is central in $\bK$ with $\hat A \not \in \bS$, then $\langle \bS, \hat A \rangle = \moddash{\bK^c}$.
\ele

\begin{proof}
    Consider the kernel $F$ of the coevaluation map corresponding to $A$ in $\bK^c$:
    \[
    0 \to F \to \hat \unit \to \hat A \otimes \hat A^*.
    \]
    Since the map $A \to A \otimes A^* \otimes A$ is the inclusion of a direct summand, one has that $\hat A \to \hat A \otimes \hat A^* \otimes \hat A$ is injective, hence $F \otimes \hat A \cong 0$. Since $\hat A$ tensor-commutes with all objects of $\bK^c$, we have $F \otimes \hat B \otimes \hat A \cong 0$ for all objects $B \in \bK^c$. Since each $G \in \moddash{\bK^c}$ is the quotient of some $\hat B$, for $B \in \bK^c$, and since the tensor product in $\moddash{\bK^c}$ is right-exact, we have $F \otimes \moddash{\bK^c} \otimes \hat A = \{0\}.$ Since $\bS$ is a Serre prime that does not contain $\hat A$, we have that $F \in \bS$. But since $\bS' = \langle \hat A, \bS\rangle$ is a Serre subcategory containing $F$ and $\hat A$, it contains $\hat \unit$. Therefore, $\bS' = \moddash{\bK^c}$. 
\end{proof}

Under certain conditions, the definition of homological primes coincides 
with the maximal Serre ideals. 

\bth{thick-max} If $\bK$ satisfies the Central Generation Hypothesis from the introduction, then the collection of homological primes of $\bK$ coincides with the collection of maximal Serre ideals in $\moddash{\bK^c}$.
\eth

\begin{proof}
Suppose all thick ideals of $\bK^c$ are centrally generated, and let $\bS$ be a homological prime in $\moddash{\bK^c}$. Set $\bP := \widetilde{\phi}(\bS)$. Suppose $F \in \moddash{\bK^c}$ and $F \not \in \bS$. Then $\langle F, \bS \rangle$ contains $\hat A$ for some $A \in \bK^c$ with $A \not \in \bP$. Since all thick ideals of $\bK$ are centrally generated, there exists $A' \in \bK^c$ with $A'$ central and $\langle A \rangle = \langle A' \rangle$. It follows that since $\hat A \in \langle F, \bS \rangle$, we also have $\hat A' \in \langle F, \bS \rangle$. By \leref{homol-prime-max}, $\langle \hat A', \bS \rangle = \moddash{\bK^c}$, and hence $\langle F, \bS \rangle = \moddash{\bK^c}$, i.e., $\bS$ is a maximal Serre ideal in $\moddash{\bK^c}$. 

For the other direction, if $\bS$ is a maximal Serre ideal of $\moddash{\bK^c}$, by \leref{inv-im}, $\widetilde{\phi}(\bS):=\bP$ is a prime ideal in $\bK$. It is certainly clear that since $\bS$ is maximal, it is also maximal in the fiber of $\widetilde{\phi}(\bS)$.
\end{proof} 

\bde{hom-spec} Define the {\em{homological spectrum}}, $\Spch \bK^c$, 
of the compact part $\bK^c$ of a rigidly-compactly generated monoidal triangulated category $\bK$ to be the collection of homological primes of $\moddash{\bK^c}$. The closed sets of the topology are defined to be arbitrary intersections of the base of closed sets
\begin{equation}
\label{Vh}
V^{\h}(A):=\{\bS \in \Spch \bK^c : \hat A \not \in \bS\}
\end{equation}
for $A \in \bK^c$. In other words, $\Spch \bK^c$ is the subsets of the Serre monoidal spectrum $\Spc(\moddash{\bK^c})$ with the induced topology. 
\ede

\section{Support maps}\label{S:supportmap}

The following proposition is an improvement of \cite[Theorem 3.14]{Vashaw-Yakimov}.

\bpr{cont-ideal-avoid-mult}
Let $\bI$ be a Serre ideal in $\moddash{\bK^c}$ and $\bM$ an nc-multiplicative subset of $\bK^c$: that is, $\bM$ satisfies the property that if $A$ and $B$ are in $\bM$, then there exists $C$ in $\bK^c$ with $A \otimes C \otimes B \in \bM$. Suppose that $\bI \cap \h(\bM) = \varnothing$. Then there exists a prime Serre ideal $\bS$ of $\moddash{\bK^c}$ maximal with respect to the property that $\bS$ contains the ideal $\bI$ and intersects $\h(\bM)$ trivially.
\epr

\begin{proof}
By Zorn's lemma, there exists a Serre ideal $\bS$ which is maximal among those containing $\bI$ and intersecting $\h(\bM)$ trivially. It remains to show that $\bS$ is prime. Suppose $F \otimes \moddash{\bK^c} \otimes G \subseteq \bS$ for some $F, G \in \moddash{\bK^c}$, and suppose for the sake of contradiction that $F$ and $G$ are not in $\bS$. By the maximality assumption on $\bS$, the Serre ideal $\langle \bS, F\rangle$ contains an object $\hat A$ for some $A \in \bM$, and likewise $\langle \bS, G\rangle$ contains an object $\hat B$ for some $B \in \bM$. One can directly verify that $\hat A \otimes \moddash{\bK^c} \otimes \hat B \subseteq \bS$. But since $\bM$ is an nc-multiplicative subset, there exists some object $C \in \bK^c$ such that $A \otimes C \otimes B \in \bM$, and hence $\hat A \otimes \hat C \otimes \hat B \in \h(\bM) \cap \bS.$ This is a contradiction. 
\end{proof}

The Serre monoidal spectrum $\Spc (\moddash{\bK^c})$ and the homological spectrum $\Spch \bK^c$ give rise to a homological support map for $\bK^c$ that has the tensor product property.

\bpr{Serre-sup} 
The map $\widetilde{V}^{\h} : \bK^c \to \Spc (\moddash{\bK^c})$ given by \eqref{wtVh} is a support map satisfying the tensor product property. In particular, the map $V^{\h} : \bK^c \to \Spch \bK^c$ given by \eqref{Vh} is a support map satisfying the tensor product property.
\epr

\begin{proof} \prref{cont-ideal-avoid-mult} implies that $\widetilde{V}^{\h}$ satisfies condition (i) in \deref{sup}. 
Condition (ii) can be directly verified. Condition (iii) 
holds because every Serre ideal of $\moddash{\bK^c}$ is 
invariant under the suspension functor \eqref{Susp}. 
$\widetilde{V}^{\h}$ satisfies condition (iv) because
the Yoneda embedding is a homological functor. 

Finally, one need to verify that $\widetilde{V}^{\h} : \bK^c \to \Spch \bK^c$ satisfies the tensor product property. For $A, B \in \bK^c$ and $\bS \in \Spc(\moddash{\bK^c})$, 
\begin{align*}
&\bS \notin \bigcup_{C \in \bK^c} \widetilde{V}^{\h} (A   \otimes C \otimes B) \; \; \Leftrightarrow \; \; 
\hat{A} \otimes \hat{C} \otimes \hat{B} \in \bS,\ \forall C \in \bK^c \; \; \Leftrightarrow 
\\ &\hat{A} \otimes \moddash{\bK^c} \otimes \hat{B} \subseteq \bS \; \; \Leftrightarrow \; \; \hat{A} \in \bS \; \mbox{or} \; \hat{B} \in \bS \; \; \Leftrightarrow \; \; 
\bS \notin \sigma (A) \cap \sigma(B),
\end{align*}
where we used that the Yoneda embedding is a monoidal functor 
(\leref{yoneda-tensor}(b)) and the fact that every object in $\moddash{\bK^c}$ is a quotient of an object in the image of the Yoneda functor. 

The map $V^{\h} : \bK^c \to \Spch \bK^c$ is a restriction of the support map $\widetilde{V}^{\h} : \bK^c \to \Spc (\moddash{\bK^c})$,
\begin{equation}
    \label{Vrestr}
V^{\h}(A) = \widetilde{V}^{\h}(A) \cap \Spch \bK^c,
\end{equation}
and by \reref{restr}, it has the stated properties. 
\end{proof}

The previous result permits us to state the following definition. 

\bde{homological-sup} We will call $V^{\h} : \bK^c \to \Spch \bK^c$ the {\em{homological support map}} 
of $\bK^c$. 
\ede

The remainder of this section will be devoted to formulating the Nerves of Steel Condition. 

\ble{ideals-inj}
The map from thick ideals of $\bK^c$ to Serre ideals of $\moddash{\bK^c}$ defined by sending a thick ideal $\bI$ to the Serre closure $\langle \h(\bI) \rangle$ of $\h(\bI)$ in $\moddash{\bK^c}$ is an injective map. In other words, if $A \not \in \bI$ in $\bK$, then $\hat A \not \in \langle \h(\bI) \rangle$ in $\moddash{\bK^c}$. 
\ele

\begin{proof}
%Since $\bI$ is a thick ideal of $\bK^c$, the Verdier quotient $\bK^c/ \bI$ is a monoidal triangulated category. 
By \leref{yoneda-tensor}(d), $\moddash{(\bK^c/\bI)} \cong \moddash{\bK^c}/ \langle \h(\bI) \rangle$. If $A \not \in \bI$, then the image of $\hat A$ in $\moddash{\bK^c/\bI}$ is nonzero, and so, $\hat A \not \in \langle \h(\bI)\rangle$. 

%We now show that $\moddash{\bK/\bI}$ and $\moddash{\bK^c}/\langle \h(\bI) \rangle$ coincide. By pulling back along the quotient $\bK \to \bK/\bI$, one obtains a functor $\moddash{\bK/\bI} \to \moddash{\bK^c}$.
\end{proof}

The restriction
\[
\phi := \widetilde{\phi}|_{\Spch \bK^c} : \Spch \bK^c \to \Spc \bK^c
\]
will be called the homological comparison map, recall \eqref{eq:phi}.

By \prref{Serre-sup}, the map $\widetilde{V}^{\h} : \bK^c \to \Spc (\moddash{\bK^c})$ and the homological support map $V^{\h} : \bK^c \to \Spc (\moddash{\bK^c})$ are support maps (in the sense of \deref{sup}) that have the tensor product property. Therefore, they should be linked to the Balmer support as in \prref{univ}. The following proposition tells us that this is done via the maps $\widetilde{\phi}$ and $\phi$, respectively.

\bpr{surj}
\hfill
\begin{enumerate}
\item[(a)] For all $A \in \bK^c$, 
\begin{equation}
\label{wtphiV}
\widetilde{\phi}^{-1}(V^B(A)) = \widetilde{V}^{\h}(A),
\end{equation}
and as a consequence,
\begin{equation}
\label{phiV}
\phi^{-1}(V^B(A)) = V^{\h}(A),
\end{equation}
recall the Balmer support \eqref{Balmer-sup}.
\item[(b)] The map $\widetilde{\phi}: \Spc(\moddash{\bK^c}) \to \Spc \bK$ is a continuous map, and as a consequence, so is the map $\phi: \Spc^{\h}\bK \to \Spc \bK$. 
\item[(c)] The map $\phi: \Spc^{\h}\bK \to \Spc \bK$ surjective. 
\end{enumerate}
\epr

\begin{proof} First we prove \eqref{wtphiV}. For all $A \in \bK^c$, 
\begin{align*}
\widetilde{\phi}^{-1}(V^B(A))&=\{ \bS \in \Spc(\moddash{\bK^c}) : A \not \in \phi(\bS)\}\\
&=\{ \bS \in \Spch(\moddash{\bK^c}) : \hat A \not \in \bS\}\\
&=\widetilde{V}^{\h}(A). 
\end{align*}
Equation~ \eqref{phiV} follows from \eqref{Vrestr} and \eqref{wtphiV}. This proves part (a). 

Part (b) follows from part (a) since the topologies of 
$\Spc( \moddash{\bK^c})$ and $\Spch \bK^c$ are given in terms 
of bases consisting of $\widetilde{V}^{\h}(A)$ and $V^{\h}(A)$ for $A \in \bK^c$. 

For part (c), let $\bP \in \Spc \bK$. To construct a homological prime $\bS$ with $\phi(\bS)=\bP$, we employ \prref{cont-ideal-avoid-mult}, taking $\bI = \langle \h(\bP) \rangle$ and $\bM=\bK \backslash \bP.$ The fact that $\bI$ and $\h(\bM)$ are disjoint follows from \leref{ideals-inj}.
\end{proof}

The developed framework for homological spectra of (the compact parts) of rigidly-compactly generated monoidal categories leads us to a point where we can formulate an extension of 
Balmer's Nerves of Steel Condition \cite{Balmer2020c}.

\bde{nerves} [{\em{Nerves of Steel Condition}}] Given a rigidly-compactly generated monoidal category $\bK$, the map 
\[
\phi : \Spch \bK^c \to \Spc \bK^c
\]
is a homeomorphism.
\ede
The remaining part of this paper is devoted to developing methods towards proving when the Nerves of Steel Condition holds. At times, we will refer to this condition as simply {\em{Nerves of Steel}}.

\bre{nerves-bijection}
Note that by \prref{surj} (b) and (c), a subset $S$ of $\Spc \bK^c$ is closed if and only if $\phi^{-1}(S)$ is closed in $\Spch \bK^c$, and so the homological comparison map $\phi$ is a quotient map, hence a homeomorphism if and only if it is bijective. Since it is already proven to be surjective, to prove Nerves of Steel it suffices to show that $\phi$ is injective. 
\ere

%%%%%%%%%%%%%%%%%%%%%%%%%%5
\section{Pure-injectives}
\label{6}
%%%%%%%%%%%%%%%%%%%%%%%%%%%%%%
Recall that in the category $\bK$, an object $A$ is called {\it pure-injective} if every morphism $f: A \to B$ such that $\hat f: \hat A \to \hat B$ is injective splits in $\bK$. The critical observation which emphasizes the importance of pure-injectives in the study of compactly generated triangulated categories is the following result of Krause \cite[Theorem 1.8, Corollary 1.9]{Krause2000} (which does not require any monoidal or tensor structure on $\bK$):

\bth{Krause-pureinj}
Let $\bK$ be a compactly generated triangulated category. Then the restricted Yoneda functor $\h$ gives a bijection between pure-injective objects of $\bK$ and injective objects of $\Moddash{\bK^c}$. Additionally, if $B$ is a pure-injective of $\bK$, there is an isomorphism
\[
\Hom_{\bK}(A,B) \cong \Hom_{\Moddash{\bK^c}} (\hat A, \hat B)
\]
for any $A \in \bK$.
\eth

For the remainder of this section, fix a rigidly-compactly generated monoidal triangulated category $\bK$, and a Serre ideal $\bS \subseteq \moddash{\bK^c}$. We will denote $\bA:=\Moddash{\bK^c}$, so that $\bA^{\fp}\cong \moddash{\bK^c}$, \eqref{eq:modMod}. Set  
\[
\overrightarrow{\bS} :=\{ F \in \bA : \text{ if } f: G \to F \text{ with } G \in \bA^{\fp}, \text{ then } f \text{ factors through }\bS\}.
\]
By \cite[Lemma 4.1]{Crawley-Boevey1994}, the subcategory $\overrightarrow{\bS}$ is the unique minimal localizing subcategory containing $\bS$, 
recall \deref{subcatsMod}(c). Since $\bS$ is a Serre ideal, $\overrightarrow{\bS}$ is also a Serre ideal of 
$\Moddash{\bK^c}$, since every object of $\overrightarrow{\bS}$ is a colimit of objects in $\bS$, and the tensor product on $\bA$ commutes with colimits.

We denote $\overline{\bA}:=\bA/\overrightarrow{\bS}$, the Gabriel--Serre quotient of $\bA$ by $\overrightarrow{\bS}$. The quotient functor
\[
Q: \bA \to \overline{\bA}
\]
admits a right adjoint
\[
R: \overline{\bA} \to \bA
\]
which is left exact, preserves injectives, and satisfies $QR \cong \Id_{\overline{\bA}}$ by \cite[Theorem 2.8]{Krause1997}. 

\bde{homol-res}
When $\bS$ is a Serre prime, the functor $Q \circ \h : \bK \to \overline{\bA}$ 
(and by abuse of notation, the abelian category $\overline{\bA}$)
will be called the {\em{noncommutative homological residue field}} associated to $\bS$. 
\ede

\bnota{bar-from-k}
We will denote the image of $F \in \bA$ under $Q$ by $\overline{F}$, and likewise denote the image of $A \in \bK$ under $Q \circ \h$ by $\overline{A}$. 
\enota

The quotient category $\overline{\bA}$ inherits a monoidal product from $\bA$ such that the quotient functor $Q$ is monoidal, by a similar argument to \cite[Appendix A]{BKS2019}:

\bpr{quot-mon}
Let $\bA$ and $\overline{\bA}$ as above. Then $\overline{\bA}$ is a monoidal category, where we define $\overline{F} \otimes \overline{G} := \overline{F \otimes G}$ for all $F, G \in \bA$.
\epr

\begin{proof}
    To prove that this monoidal product is well-defined, we must verify that if 
    \[
    0 \to F_1 \to F_2 \to F_3 \to 0
    \]
    is a short exact sequence in $\bA$ with $F_1 \in \overrightarrow{\bS}$ then $\overline{F_2 \otimes G} \cong \overline{F_3 \otimes G}$, and similarly for the roles of $1$ and $3$ reversed; and then likewise by tensoring $G$ on the right. The right and left cases are symmetric, but the cases where $1$ and $3$ are reversed are not, since the tensor product on $\bA$ is right exact, but not left exact. We have an exact sequence
    \[
    \Tor_1(F_3, G) \to F_1 \otimes G \to F_2 \otimes G \to F_3 \otimes G \to 0.
    \]
    On one hand, if $F_1 \in \overrightarrow{\bS}$, then $F_1 \otimes G$ is in $\overrightarrow{S}$ by its tensor ideal property, hence $\overline{F_2 \otimes G} \cong \overline{F_3 \otimes G}$. On the other hand, assume instead that $F_3 \in \overrightarrow{\bS}$. We must show that $\Tor_1(F_3,G) \in \overrightarrow{\bS}$. Indeed, since all objects in $\bA$ are colimits of objects in $\h(\bK^c)$ (cf. \eqref{eq:colimM}), we can find objects $\{A_i\}_{i \in I}$ and $\{B_j\}_{j\in J}$ of $\bK^c$ with an exact sequence of the form 
    \[
    \coprod_{i \in I} \widehat{A}_i \to \coprod_{j \in J} \widehat{B}_j \to G \to 0.
    \]
    Then $\Tor_1(F_3,G)$ is computed via the cohomology of the complex
    \[
    F_3 \otimes \coprod_{i \in I} \widehat{A}_i \to F_3 \otimes 
    \coprod_{j \in J} \widehat{B}_j \to 0,
    \]
    all the terms of which are in $\overrightarrow{\bS}$. Since $\overrightarrow{\bS}$ is a Serre subcategory of $\Moddash{\bK^c}$, $\Tor_1(F_3,G) \in \overrightarrow{\bS}$, which completes the proof of the proposition.
\end{proof}

The next theorem will be important for the development of theory because it shows that one can realize injective objects of $\overline{\bA}$ as images of pure injectives. 

\bth{pure-injs-in-abar}
Every injective object of the noncommutative homological residue field $\overline{\bA}$ is equal to $\overline{E}$ for some pure-injective $E \in \bK$, and there is an isomorphism
\[
\Hom_{\bK}(A,E) \cong \Hom_{\overline{\bA}}(\overline{A}, \overline{E})
\]
for all $A \in \bK$.
\eth

\begin{proof}
The proof proceeds similarly to \cite[Corollary 2.18]{BKS2019}; we include the details for the ease of the reader. Given an injective object of $\overline{\bA}$, we can form the pure-injective $E$ by first applying $R$ and then applying \thref{Krause-pureinj}. For the isomorphism of hom sets, compute
\begin{align*}
    \Hom_{\bK}(A,E) & \cong \Hom_{\bA} (\hat A, \hat E)\\
    &\cong \Hom_{\bA}(\hat A, R(\overline{E}))\\
    &\cong \Hom_{\overline{\bA}} (\overline{A}, \overline{E}).
\end{align*}
The first isomorphism is by \thref{Krause-pureinj}, the second by the definition of $\hat E$, and the third by the adjunction between $Q$ and $R$. 
\end{proof}

Since $\bA$ is a Grothendieck category, the noncommutative homological residue field $\overline{\bA}$ is also a Grothendieck category, and each such category has injective envelopes \cite[Theorem 6.25]{Freyd}.

\bnota{pure-inj-s}
We will denote the injective envelope of $\overline{\unit}$ in $\overline{\bA}$ by $\overline{E_{\bS}},$ and the unique pure-injective which corresponds to $\overline{E_{\bS}}$ in $\bK$ by \thref{pure-injs-in-abar} will be denoted $E_{\bS}$.
\enota

In the remainder of this section, since $\bS$ is fixed, we will denote \[
E:=E_{\bS}.
\]
Our main goal will be showing that the localizing category $\overrightarrow{\bS}$ is equal to the kernel of the functor $\hat{E} \otimes -$ in $\Moddash{\bK^c}$. 

We begin by proving a noncommutative version of \cite[Proposition 3.3]{BKS2019}. Note that by \thref{pure-injs-in-abar} we have a distinguished triangle
\[
A \xrightarrow{f} \unit \xrightarrow{g} E \to \Sigma A
\]
such that $\overline{g}: \overline{\unit} \to \overline{E}$ is the injective envelope for $\overline{\unit}$ in $\overline{\bA}$. 

\ble{f-inj-properties-in-k}
Let $f$ and $g$ be as above. Then
\begin{enumerate}
    \item[(a)] if $h: B \to \unit$ in $\bK$ satisfies $\overline{h}=0$, then there exists $i: B \to A$ such that $h=fi$;
    \item[(b)] $f \otimes \id_E$ and $\id_E \otimes f$ are both 0;
    \item[(c)] if $h: B \to C$ is a morphism in $\bK$ such that $C \in \bK^c$ and $\overline{h}=0$, then $h \otimes \id_E$ and $\id_E \otimes h$ are both 0.
\end{enumerate}
\ele

\begin{proof} (a) By the extension axiom for triangulated categories, it is enough to show that $gh=0$, via the morphism of triangles
\begin{center}
    \begin{tikzcd}
A \arrow[r, "f"]                                    & \unit \arrow[r, "g"]       & E \arrow[r]                & \Sigma A \\
B \arrow[r, "\id_B"] \arrow[u, "\exists i", dashed] & B \arrow[u, "h"] \arrow[r] & 0 \arrow[u, "0"] \arrow[r] & \Sigma B
\end{tikzcd} 
\end{center}
    Now note that $\overline{gh} = \overline{g} \overline{h} =0$, since $\overline{h}$ is 0 by assumption; but by the isomorphism $\Hom_{\bK}(B,E) \cong \Hom_{\overline{\bA}}(\overline{B}, \overline{E})$ from \thref{pure-injs-in-abar}, this implies that $gh=0$, which proves (a). 

    (b) Since $\h$ is a homological functor, and since the quotient functor $\bA \to \overline{\bA}$ is exact, the sequence
    \[
    \overline{A} \xrightarrow{\overline{f}} \overline{\unit} \xrightarrow{\overline{g}} \overline{E}
    \]
    is exact in $\overline{\bA}$. Since $\overline{g}$, by assumption, is the map of the unit into its injective envelope, it is an injective map and $\overline{f}$ is 0. This applies immediately that $\overline{f \otimes \id_E},$ which is identified with $\overline{f} \otimes \id_{\overline{E}}$, is 0 as a morphism $ \overline{A} \otimes \overline{E} \to \overline{E}$. Again using the isomorphism from \thref{pure-injs-in-abar}, this implies $f \otimes \id_E =0$. The argument that $\id_E \otimes f =0$ follows in the same way.

    (c) One has 
    \[
    h':= C^* \otimes B \xrightarrow{\id_{C^*} \otimes h} C^* \otimes C \xrightarrow{\ev_C} \unit
    \]
    is a morphism whose target is $\unit$, and satisfies $\overline{h'}=0$. By (a) $h'$ factors through $f$, and hence by (b), $h' \otimes \id_E$ and $\id_E \otimes h'$ are 0. Now just note that by the defining equations for duals, $h$ is equal to
    \[
    B \xrightarrow{\coev_C \otimes \id_B} C \otimes C^* \otimes B \xrightarrow{\id_C \otimes h'} C.
    \]
    Hence, $h \otimes \id_E$ factors through $\id_C \otimes h' \otimes \id_E=0$, and so $h \otimes \id_E=0$. 
    
    On the other hand, define \[     h'':=B \otimes {^* C} \xrightarrow{h \otimes \id_{{^* C}}} C \otimes {^*C} \xrightarrow{\ev_{{^*C}} }\unit      \] 
    using the right dual ${^*C}$ of $C$. We have that $\overline{h''}=0$, hence $h''$ factors through $f$, and so $h'' \otimes \id_E$ and $\id_E \otimes h''=0$. Again using the duality conditions, we can write $h$ as
    \[
    B \xrightarrow{\id_B \otimes \coev_{{^*C}}} B \otimes {^* C} \otimes C \xrightarrow{h'' \otimes \id_C} C.
    \] Since $\id_E \otimes h''=0$, it follows that $\id_E \otimes h=0$, and (c) as follows. 
\end{proof}

\bth{S-is-KerE} One has equalities of the following localizing subcategories in $\bA$: 
\[
\ker(\hat E \otimes -) = \overrightarrow{\bS} = \ker(- \otimes \hat E).
\]
\eth

\begin{proof}
    The proof is similar to \cite[Theorem 3.5]{BKS2019}, with a few necessary modifications. First, suppose that $F \in \bA$ satisfies $\hat E \otimes F = 0$. Recall that since every object of $\bA$ is a colimit of objects from $\h(\bK^c)$, we can take objects $\{A_i\}_{i \in I}$ and $\{B_j\}_{j\in J}$ of $\bK^c$ together with maps forming an exact sequence \[
    \coprod_{i \in I} \hat A_i \to \coprod_{j \in J} \hat B_j \to F \to 0.
    \]
    Denote the coproduct of the $A_i$'s by $A$ and the coproduct of the $B_j$'s by $B$, so that we have maps $f: A \to B$, $g: \hat B \to F$ such that 
    \[
    \hat A \xrightarrow{\hat{f}} \hat B \xrightarrow{g} F \to 0
    \]
    is exact in $\bA$. 
    
    Since 
    \[
    \hat A \xrightarrow{\hat f} \hat B \xrightarrow{\hat h} \hat C
    \]
    is also exact, for $C=\cone(f)$, it follows that $F$ is the image of $\hat h$. Since tensor product with $\hat E$ is exact, the image of $\id_{\hat E} \otimes \hat h =  \hat E \otimes F \cong 0$, and so $\id_{\hat E} \otimes \hat h = 0$. In particular, the image in $\overline{\bA}$ also vanishes: $\id_{\overline{E}} \otimes \overline{h}=0$. 
    
    Recall that $\overline{E}$ is the injective envelope of $\overline{\unit}$ in $\overline{\bA}$, say via 
    \[
    i: \overline{\unit} \to \overline{E}.
    \]
    There exists a commutative diagram
    \[
    \begin{tikzcd}
\overline{\unit} \otimes \overline{B} \arrow[d, "i \otimes \id_{\overline{B}}"'] \arrow[r, "\id_{\overline{\unit}} \otimes \overline{h}"] & \overline{\unit} \otimes \overline{C} \arrow[d, "i \otimes \id_{\overline{C}}"] \\
\overline{E} \otimes \overline{B} \arrow[r, "\id_{\overline{E}} \otimes \overline{h} = 0"']                                               & \overline{E} \otimes \overline{C}                                              
\end{tikzcd}
    \]
where the map $i \otimes \id_{\overline{C}}$ is injective, since $\overline{C}$ is flat and $i$ is injective; this immediately implies that $\id_{\overline{\unit}} \otimes \overline{h} =0$. It follows that $\overline{h} = 0$, and since $\overline{F}$ is the image of $\overline{h}$, 
$\overline{F} \cong 0$, that is, $F \in \overrightarrow{\bS}$. This shows that $\ker(\hat E \otimes -) \subseteq \overrightarrow{\bS}$; the argument that $\ker(- \otimes \hat E) \subseteq \overrightarrow{\bS}$ is symmetric.

Next, we show $\overrightarrow{\bS} \subseteq \ker(\hat E \otimes -) \cap \ker(- \otimes \hat E).$ Note that $\hat E$ is flat, and so $\ker(\hat E \otimes -)$ is a localizing subcategory of $\bA$. Since $\overrightarrow{\bS}$ is generated as a localizing subcategory by $\bS$, it follows that $\overrightarrow{\bS} \subseteq \ker(\hat E \otimes -)$ if and only if $\bS \subseteq \ker(\hat E \otimes -)$. Let $F \in \bS.$ As before, we can find a morphism $f: A \to B$ in $\bK^c$ such that $F$ is the image of $\hat{f}$. But since the quotient functor $\bA \to \overline{\bA}$ is exact, the image of $\overline{f}$ is $\overline{F}$, and $\overline{F} \cong 0$ in $\overline{\bA}$ (using the assumption that $F \in \bS$). It follows that $\overline{f}: \overline{A} \to \overline{B}$ is equal to $0$. By \leref{f-inj-properties-in-k}(c), using the fact that $B$ is compact, we can conclude that $f \otimes \id_E=0=\id_E \otimes f$. Since $F$ is the image of $\hat{f}$, it immediately follows that $F \otimes \hat E\cong 0 \cong \hat E \otimes F$, showing that $F \in \ker(\hat E \otimes -) \cap \ker(- \otimes \hat E)$.

We have now shown that 
\[
\ker(\hat E \otimes -) \cup \ker(- \otimes \hat E) \subseteq \overrightarrow{\bS} \subseteq \ker(\hat E \otimes -) \cap \ker(- \otimes \hat E).
\]
The claim of the theorem follows.
\end{proof}

%%%%%%%%%%%%%%%%%%%%%%%%%%%%%%%%%%%%%%%%%%%%%%%%%%%%%%%
\section{Stratification vs Nerves of Steel}
\label{7}
%%%%%%%%%%%%%%%%%%%%%%%%%%%%%%%%%%%%%%%%%%%%%%%%%%%5%%

In this section, we connect the Nerves of Steel Condition to the stratification of $\bK$ via the Balmer--Favi support. Recall the construction of extended supports inspired by Benson--Iyengar--Krause \cite{BIK2008} and Balmer--Favi \cite{Balmer-Favi}, which was fully realized in the noncommutative case in \cite{CaiVashaw}. 
The {\em{extended support}} is a map which assigns to each object of $\bK$ a subset of $\Spc \bK^c$, such that if $A$ is compact the subset is the Balmer support of $A$, which satisfies typical properties of support. Moreover, 
the extended support is defined using the Rickard idempotent functors $L_{\bI}$ and $\Gamma_{\bI}$, which are triangulated functors $\bK \to \bK$, associated to a thick subcategory $\bI$ of $\bK$; see e.g.,~ \cite[Section 3]{CaiVashaw} for background on these functors. Critically, for any object $A \in \bK$, there is a distinguished triangle
\begin{equation}
\label{rickard-tri}
    \Gamma_{\bI} A \to A \to L_{\bI} A \to \Sigma \Gamma_{\bI} A.
\end{equation}
The object $\Gamma_{\bI} A$ is in $\Loc(\bI)$ and the object $L_{\bI} (A)$ is in $\Loc(\bI)^{\perp}$, that is, it admits no maps from any objects of $\Loc(\bI)$; and if 
\[
B \to A \to C \to \Sigma B
\]
is another triangle with $B \in \Loc(\bI)$ and $C \in \Loc(\bI)^{\perp}$, then there are unique isomorphisms $B \cong \Gamma_{\bI}(A)$ and $C \cong L_{\bI}(A)$. Using this universal property, one can verify that 
\[
\Gamma_{\bI}(A) \cong \Gamma_{\bI}(\unit) \otimes A\cong A \otimes \Gamma_{\bI}(\unit), \quad L_{\bI}(A) \cong L_{\bI}(\unit) \otimes A\cong A \otimes L_{\bI}(\unit). 
\]

Recall that given any specialization-closed subset $S$ of $\Spc \bK^c$, 
one has a thick ideal
\[
\Theta^B(S):=\{ A \in \bK^c : V^B(A) \subseteq S\}.
\]
For ease of notation, denote
\[
\Gamma_{S}:=\Gamma_{\Theta^B(S)}, \quad L_{S} := L_{\Theta^B(S)}
\]
for any specialization-closed subset $S$ of $\Spc \bK^c$. Now define for any $\bP$ in $\Spc \bK^c$ the functor
\[
\Gamma_{\bP} :=\Gamma_{\mc{V}(\bP)} L_{\mc{Z}(\bP)} \cong L_{\mc{Z}(\bP)} \Gamma_{\mc{V}(\bP)},
\]
where $\mc{V}(\bP)$ and $\mc{Z}(\bP)$ are the subsets
\[
\mc{V}(\bP) := \overline{\{ \bP\}}, \quad \mc{Z}(\bP) := \{ \bQ \in \Spc \bK^c : \bP \not \in \mc{V}(\bQ)\}.
\]
The extended support, as defined in \cite{CaiVashaw} and which we will refer to throughout the remainder of the paper as the {\em{Balmer--Favi support}}, is defined as
\begin{equation}
\label{BFsupport}
V^{BF}(A) := \{ \bP \in \Spc \bK^c : \Gamma_{\bP} A \not \cong 0 \}
\end{equation}
for any $A \in \bK$. Note that if $A$ is compact, the Balmer--Favi support recovers the Balmer support, that is, $V^B(A) = V^{BF}(A)$, see \cite[Proposition 3.8]{CaiVashaw}.

Let $\bL$ be a localizing ideal of $\bK$ and $S$ be a subset of $\Spc \bK^c$. We denote 
\[
\Phi^{BF}(\bL):= \bigcup_{A \in \bL} V^{BF}(A)
\]
and 
\[
\Theta^{BF}(S) :=\{A \in \bK : V^{BF}(A) \subseteq S\}.
\]

\bde{strat} The topological space $\Spc \bK^c$ {\it{stratifies}} $\bK$ if the maps $\Phi^{BF}$ and $\Theta^{BF}$ define an order-preserving bijection between localizing $\otimes$-ideals of $\bK$ and subsets of $\Spc \bK^c$.
\ede

The next theorem is an explicit description of the Balmer--Favi support of all pure-injective objects representing arbitrary Serre primes of $\moddash{\bK^c}$. Our result is new even in the symmetric case and is a far-reaching extension of \cite[Lemma 3.7]{BHS2023}.

\bth{supp-pureinj-2}
Suppose $\Spc \bK^c$ is Noetherian. Let $\bS$ be a Serre prime of $\moddash{\bK^c}$ with corresponding pure-injective 
\[
E:=E_{\bS} \in \bK,
\]
which satisfies $\overrightarrow{\bS} = \ker (\hat{E} \otimes -).$ Then the Balmer--Favi support for $E$ is equal to 
\[
V^{BF}(E)=\{ \phi(\bS') : \bS' \in \Spch(\bK^c) \; \text{ with } \; \bS' \supseteq \bS \}.
\]
In other words, the Balmer--Favi support of $E$ consists 
of the primes of $\bK^c$ living over $\bS$, recall the terminology in the introduction.
\eth

\begin{proof} Set $\bP:=\phi(\bS).$
We use throughout the proof that $\bQ$ is in $V^{BF}(E)$ if and only if $\Gamma_{\bQ} E \not \cong 0$, which is true if and only if $\widehat{\Gamma_{\bQ} \unit} \not \in \overrightarrow{\bS}$.
The proof will consist of the following steps.
\begin{enumerate}
    \item $\bP$ is in the Balmer--Favi support of $E$;
        \item If there exists $\bS' \in \Spch(\bK^c)$ with $\bS \subseteq \bS'$, then $\bQ:=\phi(\bS')$ is in the Balmer--Favi support of $E$;
    \item If $\bQ$ is in the Balmer--Favi support of $E$, then $\bQ$ contains $\bP$;
    \item If $\bQ$ is in the Balmer--Favi support of $E$, then there exists $\bS'$ in $\Spch(\bK^c)$ containing $\bS$ such that $\bQ = \phi(\bS')$. 
\end{enumerate}

For (1), note that $\bP \in V^{BF}(E)$, i.e.,~ that $\Gamma_{\bP} E \not \cong 0,$ or equivalently that $\widehat{\Gamma_{\bP} \unit} \not \in \overrightarrow{\bS}$. We first claim that $L_{\mc{Z}(\bP)} E \cong E$, that is, that $E $ is in $\Loc(\Theta^B(\mc{Z}(\bP)))^{\perp}$. Indeed, if $A$ is compact then $A \in \Theta^B(\mc{Z}(\bP))$ if and only if $V^B(A) \subseteq \mc{Z}(\bP)$ (that is, $\bP \not \in V^B(A)$, i.e.,~ $A \in \bP$). We know in this case that $A \otimes E \cong 0$, and likewise $A^* \otimes E \cong 0$; since $A$ is compact, this implies
\[
\Hom_{\bK}(A,E) \cong \Hom_{\bK}(\unit, A^* \otimes E) \cong \Hom_{\bK}(\unit,0) \cong 0.
\]
Hence, $E$ is indeed in $\Loc(\Theta^B(\mc{Z}(\bP)))^{\perp}$, and so $L_{\mc{Z}(\bP)} E \cong E$. 

To conclude that $\Gamma_{\bP} E \not \cong 0$, it remains to show that $\Gamma_{\mc{V}(\bP)} E \not \cong 0$. Recall that $\Gamma_{\mc{V}(\bP)} E \cong 0$ if and only if $E \in \Loc(\Theta^B(\mc{V}(\bP)))^{\perp}$. This will hold if and only if $E \otimes A \cong 0$ for all $A \in \Theta^B(\mc{V}(\bP))$. To contradict this, take $B$ to be a compact object with $V^{BF}(B)=\mc{V}(\bP)$; such a $B$ exists by \thref{thm-bij}, using the fact that $\Spc \bK^c$ is Noetherian. Then certainly $B \in \Theta^B(\mc{V}(P))$, and additionally $B \not \in \bP$. But since $B \not \in \bP$, one also has $\hat B \not \in \bS$, and so $E \otimes B \not \cong 0$. Therefore, $E$ is not in $\Loc(\Theta^B(\mc{V}(\bP)))^{\perp}$, and (1) follows.  

Now assume the situation of (2). Set $E'$ to be the pure-injective corresponding to $\bS'$. Suppose $\bQ$ were not in $V^{BF}(E)$. This would mean by definition that $\Gamma_{\bQ}\unit \otimes E \cong 0$, hence that $\Gamma_\bQ 
 \unit $ is in $\overrightarrow{\bS}$. However, $\Gamma_{\bQ}\unit$ would be in $\overrightarrow{\bS'}$, and hence $\Gamma_{\bQ} \unit \otimes E'$ would also be 0. This would then contradict (1) applied to $E'$, since $\bQ$ must be in the support of $E'$.

For (3), suppose that $\bQ$ is a prime ideal which does not contain $\bP$. We need only show that $\Gamma_{\mc{V}(\bQ)} \unit \in \Loc(\bP)$. But indeed, $\Gamma_{\mc{V}(\bQ)} \unit$ is in $\Loc(\Theta^B(\mc{V}(\bQ)))$, so it suffices to prove that $\Theta^B(\mc{V}(\bQ)) \subseteq \bP$. By definition, $\Theta^B(\mc{V}(\bQ))$ is the thick ideal consisting of compact objects $A$ with $V^B(A) \subseteq \mc{V}(\bQ) = \overline{\{\bQ\}}$. Since $\bP \not \subseteq \bQ$, we have $\bP \not \in \overline{\{\bQ\}}$; hence if $A \in \Theta^B(\mc{V}(\bQ))$, $\bP \not \in V^B(A)$, i.e.,~ $A \in \bP$. This completes the proof of claim (3).

Finally, for (4), suppose $\bQ \in V^{BF}(E)$ and $\bQ \not = \bP$. By (3), we know that $\bQ$ strictly contains $\bP$. Consider 
\[
\bM:=\{ \hat A \in \moddash{\bK^c} : A \in \bK^c \backslash \bQ\},
\]
and consider
\[
\bS':=\langle \bS, \hat{A} : A \in \bQ\rangle.
\]
One can readily verify that $\bM$ is an nc-multiplicative subset, using the prime condition. \prref{cont-ideal-avoid-mult} can be used to produce a Serre prime $\bS''$ which is maximal with respect to the property that $\bS''$ contains $\bS'$, and intersects $\bM$ trivially. If one can do this, it is clear by construction that $\phi(\bS'')=\bQ$, and that $\bS'' \in \Spch(\bK^c)$. 

To apply \prref{cont-ideal-avoid-mult}, we simply need to verify that $\bM$ and $\bS'$ are disjoint. Suppose for the sake of contradiction that there is an object $C \in \bK^c \backslash \bQ$, such that $\hat C \in \bS'$. Consider the Serre ideal $\overline{\bS'} = \{ \overline{F} : F \in \bS'\} = \langle \overline{A} \in \overline{\bA} : A \in \bQ \rangle$ 
of the noncommutative homological residue field
$\overline{\bA}$. 
Note that for any object $\overline{F}$ of $\overline{\bS}'$, we have $\overline{\Gamma_{\bQ} \unit} \otimes \overline{F} \cong 0$, since $\Gamma_{\bQ} \unit \otimes A \cong 0$ for all $A \in \bQ$; since $\Gamma_{\bQ} \unit $ is idempotent, $\overline{\Gamma_{\bQ}\unit}$ cannot be in the localizing category generated by $\overline{\bS}'$. In particular, this means that $\widehat{\Gamma_{\bQ}\unit}$ cannot be in $\overrightarrow{\bS'}$. But indeed, $\Gamma_{\bQ} \unit$ is in $\Loc(\langle C \rangle)$: $\Loc(\langle C \rangle)$ is a tensor ideal, $\Gamma_{\mc{V}(\bQ)} \unit$ is in $\Loc(\Theta^B(\mc{V}(\bQ)))$, and $\Theta^B(\mc{V}(\bQ)) \subseteq \langle C \rangle$ by the assumption that $\bQ \in V^B(C)$. However, if $\Gamma_{\bQ} \unit$ is in $\Loc(\langle C \rangle)$, then it follows that $\widehat{\Gamma_{\bQ} \unit } \in \overrightarrow{\bS'}$, which is a contradiction. It follows that $\bS'$ is indeed disjoint from $\bM$, and we are done. 
\end{proof}

Theorems \ref{tthick-max} and \ref{tsupp-pureinj-2} imply the following corollary. 

\bco{supp-pureinj} Assume that $\Spc \bK^c$ is Noetherian.
\begin{enumerate}
\item[(a)] If $\bS$ is a maximal Serre ideal of $\moddash{\bK^c}$ with corresponding pure-injective $E_{\bS} \in \bK$,
then the Balmer--Favi support of $E_{\bS}$ is $\{\phi(\bS)\}$. 
\item[(b)] If all Serre primes in $\Spch(\bK^c)$ are maximal, then $\bK$ satisfies the Uniformity Hypothesis.
\item[(c)] In particular, if $\bK$ satisfies the Central Generation Hypothesis, then it satisfies the Uniformity Hypothesis.
\end{enumerate}
\eco

We will need the following tensor product property for pure-injectives to state the main result of this section; this tensor product property is likely of independent interest, since the tensor product property in monoidal triangulated categories has received significant attention in recent years.

\bpr{tpp-pureinj}
Suppose $\bK$ is stratified by $\Spc \bK^c$ and that $\Spc \bK^c$ is Noetherian. Let $\bS$ and $\bS'$ be two Serre ideals of $\moddash{\bK^c}$ with corresponding pure-injectives $E$ and $E'$. Then
\[V^{BF}(E' \otimes E) = V^{BF}(E) \cap V^{BF}(E').\]
\epr

\begin{proof}
    The containment $V^{BF}(E' \otimes E) \subseteq V^{BF}(E) \cap V^{BF}(E')$ holds by \cite[Proposition 3.8]{CaiVashaw}. For the reverse containment, suppose that $\bP \in V^{BF}(E) \cap V^{BF}(E'),$ that is, $\Gamma_{\bP} \unit \otimes E \not \cong 0$ and $\Gamma_{\bP} \unit \otimes E' \not \cong 0$. We must show that $\Gamma_{\bP} \unit \otimes E' \otimes E \not \cong 0$. Recall that $\Gamma_{\bP} \unit$ tensor-commutes with all objects of $\bK$ (see \cite[Lemma 3.2]{CaiVashaw}). By stratification, using the fact that $V^{BF}(\Gamma_{\bP} \unit) = \{\bP\}$ (\cite[Lemma 4.2, Lemma 4.3]{CaiVashaw}), it follows that $\Gamma_{\bP} \unit \in \Loc(E)$. By an analogous argument to \cite[Lemma 3.1.2]{NVY1}, it follows that
    \[
    \Gamma_{\bP} \unit \otimes E' \in \Loc(E \otimes A \otimes E' : A \in \bK).
    \]
    But since $\Gamma_{\bP} \unit \otimes \Gamma_{\bP} \unit \otimes E' \cong \Gamma_{\bP} \unit \otimes E' \not \cong 0$, it follows that there must be $A \in \bK$ such that $\Gamma_{\bP} \unit \otimes E \otimes A \otimes E' \not \cong 0.$ Using the fact that $\Gamma_{\bP} \unit$ commutes with all objects and the fact that $\ker(E \otimes -)=\ker(- \otimes E)$ from \thref{S-is-KerE}, it follows that $\Gamma_{\bP} \unit \otimes E' \otimes E \not \cong 0$, so we are done.
\end{proof}

\bco{max-loc-contain}
Suppose $\bK$ is stratified by $\Spc \bK^c$ and that $\Spc \bK^c$ is Noetherian. Let $\bS \in \Spch \bK^c$ with pure-injective $E$, and let $\bS'$ be a Serre prime in $\moddash{\bK^c}$ with pure-injective $E'$. Suppose $E \in \Loc(E')$. Then $\bS' \subseteq \bS$. 
\eco

\begin{proof}
    By \prref{tpp-pureinj}, $V^{BF}(E \otimes E')= V^{BF}(E).$ Set $\bS'':=\ker(\hat E \otimes \hat E' \otimes -)^{\fp}.$ By the stratification assumption, $\phi(\bS'')=\phi(\bS)$, since $E \otimes E'$ and $E$ generate the same localizing ideal of $\bK$. Note that from their definitions
    \[
    \overrightarrow{\bS''} \supseteq \overrightarrow{\bS} \cup \overrightarrow{\bS'},
    \]
    using the fact that $\ker(\hat E \otimes -) = \ker(- \otimes \hat E)$ from \thref{S-is-KerE}. But since $\bS \in \Spch(\bK^c)$ and is therefore maximal in the fiber of $\phi(\bS)$, it then follows that $\bS = \bS'',$ and consequently that $\bS' \subseteq \bS$.
\end{proof}

\bth{strat-nos}
    Suppose $\bK$ is stratified by $\Spc \bK^c$, satisfies the Uniformity Hypothesis from the introduction, and that $\Spc \bK^c$ is Noetherian.
    Then $\bK^c$ satisfies the Nerves of Steel Condition.
\eth

\begin{proof}
    Suppose $\bS$ and $\bS'$ are two Serre ideals maximal in the fiber of $\bP \in \Spc \bK^c$. Let $E$ and $E'$ be the pure-injectives corresponding to $\bS$ and $\bS'$, respectively. By \thref{supp-pureinj-2}, the Balmer--Favi support of $E$ consists precisely of the primes which live over $\bS$; by assumption, this is equal to the set of primes that live over $\bS'$, that is, the Balmer--Favi support of $E'$. By the stratification assumption, since the Balmer--Favi supports of $E$ and $E'$ coincide, we have $\Loc(E)=\Loc(E')$. This immediately implies by \coref{max-loc-contain} that 
    \[
    \overrightarrow{\bS} = \ker( \hat E \otimes -) = \ker(\hat E' \otimes -) = \overrightarrow{\bS'}.
    \]
    Hence, $\bP$ has a unique homological prime in its fiber, and so $\phi$ is a bijection.
\end{proof}

The following corollary generalizes \cite[Theorem 4.17]{BHS2023} from the symmetric case to the case where $\bK^c$ satisfies the central generation condition. It follows directly from 
\thref{strat-nos} and \coref{supp-pureinj}.

\bco{nos-central} Suppose $\bK$ is stratified by $\bK^c$, satisfies the Central Generation Hypothesis and
$\Spc \bK^c$ is Noetherian. Then the Nerves of Steel Condition holds for $\bK$.
\eco

%\begin{proof} The corollary follows from \thref{strat-nos} and \leref{thick-max}, since if $\bP \in \Spc \bK^c$ and $\bS \in \phi^{-1}(\bP)$, then since $\bS$ is maximal $\bP$ is the only prime living over $\bS$.
%\end{proof}
%%%%%%%%%%%%%%%%%%%%%%%%%%%%%%%%
\section{Stratification and Nerves of Steel for crossed products}
\label{8}
%%%%%%%%%%%%%%%%%%%%%%%%%%%%
Throughout this section, assume that $\Spc \bK^c$ is Noetherian, and let $G$ be a group acting on $\bK$. That is, for each $g \in G$, there exists a monoidal triangulated autoequivalence $T_g: \bK \to \bK$ with compatible natural isomorphisms $T_g T_h \cong T_{gh}$. Given a thick ideal $\bI$ of $\bK^c$ and $g \in G$, denote 
\[
T_g \bI :=\{ T_g A : A \in \bI\},
\]
which is another thick ideal of $\bK^c$. This induces an action of $G$ on $\Spc \bK^c$; for $g \in G$ and $\bP \in \Spc \bK^c$, write $g.\bP$ for the action of $g$ on $\bP$. 

Recall the crossed product category $\bK \rtimes G$ \cite[Section 4.15]{EGNO}. In this category, objects are formal direct sums of objects of the form $A \boxtimes g$ for $g \in G$ and $A \in \bK$, and morphisms are defined by setting
\[
\Hom_{\bK \rtimes G} (A \boxtimes g, B \boxtimes h) := \begin{cases} \Hom_{\bK} (A,B) & g = h, \\ 0 & g \not = h. \end{cases}
\]
The crossed product category $\bK \rtimes G$ is compactly-generated with compact part $\bK^c \rtimes G$. The monoidal product on $\bK \rtimes G$ is defined via
\[
(A \boxtimes g) \otimes (B \boxtimes h) :=(A \otimes T_g B) \boxtimes gh
\]
for $A, B \in \bK$ and $g,h \in G$. 

Crossed product categories have been important in recent years in constructing counterexamples to the tensor product property for support varieties \cite{Benson-Witherspoon,BPW2024} and counterexamples to classifying thick ideals in monoidal triangular geometry \cite{HuangVashaw}.

\bex{plavnik-witherspoon}
If $\bK= \underline{\Mod}(H)$ is the stable module category a finite-dimensional Hopf algebra $H$ over a field $k$ which admits an action by a finite group $G$, then $\bK \rtimes G$ is equivalent to the stable category of modules for the Plavnik--Witherspoon \cite{Plavnik-Witherspoon} smash coproduct $(H^* \# k G)^*$. We will explore this example in more detail in \exref{PW-Hopf} below.
\eex

We will use freely the description of $\Spc (\bK^c \rtimes G)$ given in \cite[Theorem 1.1]{HuangVashaw}. Namely, $\Spc (\bK^c \rtimes G)$ is homeomorphic to the topological space where points are defined to be $G$-orbits of points of $\Spc \bK^c$ under the induced action, and closed subsets are generated by the orbits of all closed subsets of $\Spc \bK^c$. We will write $G.\bP$ to denote the $G$-orbit of $\bP$ for $\bP \in \Spc \bK^c$. Under this description of $\Spc (\bK^c \rtimes G)$, the Balmer support of $A \boxtimes g$ is equal to the $G$-orbit of the Balmer support of $A$ in $\Spc \bK^c$, for all $A \in \bK^c$ and $g \in G$ \cite[Proposition 7.3]{HuangVashaw}. By \cite[Theorem 7.6]{CaiVashaw}, the Balmer--Favi support has the same description; that is, if $A \in \bK$ (no longer necessarily compact), then the Balmer--Favi support of $A \boxtimes g \in \Spc(\bK^c \rtimes G)$ is equal to the collection of points \[
\{G.\bP : \bP \in V^B (A) \subseteq \Spc \bK^c\}.
\]
We proceed to describe the homological spectrum of $\bK^c \rtimes G$. The following lemma describes the compatibility between the Rickard idempotent functors and the $G$-action.

\ble{tg-loc}
Let $G$ be a group acting on $\bK$. Then $\Loc(T_g \bI) = T_g \Loc(\bI)$ for any thick ideal $\bI$ of $\bK^c$. As a consequence, 
$$\Gamma_{T_g \bI} T_g A \cong T_g \Gamma_{\bI}(A)$$ 
and 
$$L_{T_g \bI} T_g A \cong  T_g L_{\bI} A$$ 
for any object $A$. 
\ele

\begin{proof}
    It is clear that $T_g \Loc(\bI)$ is a localizing subcategory of $\bK$ that contains $T_g \bI$, so the inclusion
    \[
    \Loc(T_g \bI) \subseteq T_g \Loc(\bI)
    \]
    is automatic. On the other hand, note that by the same argument
    \[
    \Loc(\bI) \subseteq T_{g^{-1}} \Loc( T_{g} \bI).
    \]
    By applying $T_g$ to this identity, we obtain $T_g \Loc(\bI) \subseteq \Loc(T_g \bI)$. The consequence follows from the universal property of the Rickard triangle (\ref{rickard-tri}), by applying $T_g$ to
    \[
    \Gamma_{\bI} A \to A \to L_{\bI} A.
    \]
\end{proof}

\ble{gamma-prime-tg}
Let $G$ be a group acting on $\bK$. Then for any $\bP \in \Spc \bK^c$, 
$$T_g \Gamma_{\bP} A \cong \Gamma_{g.\bP} T_g A.$$.
\ele

\begin{proof}
    By definition, $\Gamma_{\bP} = \Gamma_{\mc{V}(\bP)} L_{\mc{Z}(\bP)}$. One can directly verify that $g.\mc{V}(\bP)=\mc{V}(g.\bP)$ and $g.\mc{Z}(\bP)=\mc{Z}(g.\bP)$. Then by \leref{tg-loc}, 
    \begin{align*}
        T_g \Gamma_{\bP} A &\cong T_g \Gamma_{\mc{V}(\bP)} L_{\mc{Z}(\bP)} A \\
        & \cong \Gamma_{\mc{V}(g.\bP)} T_g L_{\mc{Z}(\bP)} A\\
        & \cong \Gamma_{\mc{V}(g.\bP)} L_{\mc{Z}(g.\bP)} T_g A\\
        & \cong \Gamma_{g.\bP} T_g A.
    \end{align*}
\end{proof}

We can now provide information about stratification of crossed product categories. 
 
\bth{strat-crossprod}
Suppose that $\Spc \bK^c$ is Noetherian, and let $G$ be a group acting on $\bK$. Suppose $\Spc \bK^c$ stratifies $\bK$. Then $\Spc (\bK^c \rtimes G)$ stratifies $\bK \rtimes G$. 
\eth

\begin{proof}
    First, note that inflation and restriction provide a bijection between the localizing ideals of $\bK \rtimes G$ and the $G$-invariant localizing ideals of $\bK$, by the same argument as in \cite[Proposition 5.4]{HuangVashaw}. Using the description of $\Spc(\bK^c \rtimes G)$ from \cite{HuangVashaw} and the description of the Balmer--Favi support on $\bK \rtimes G$ from \cite{CaiVashaw}, both recalled above, to prove the theorem it suffices to show that if $\bL$ is a localizing ideal of $\bK$, then $\bL$ is $G$-invariant if and only if $\Phi^{BF}(\bL)$ is $G$-invariant. 
    
    We now check that for any $A \in \bK$ and $g \in G$, we have
    \begin{align*}
        g . V^{BF}(A) &= \{ g.\bP : \bP \in V^{BF}(A)\}\\
        &= \{ g .\bP : \Gamma_{\bP} A \not \cong 0\} \\
        &= \{\bP : \Gamma_{g^{-1}.\bP} A \not \cong 0\}\\
        &= \{ \bP : T_{g^{-1}} \Gamma_{\bP} T_{g} A \not \cong 0\}\\
        & = \{ \bP : \Gamma_{\bP} T_g A \not \cong 0\}\\
        &= V^{BF}(T_g A),
    \end{align*}
    where the fourth equality is by \leref{gamma-prime-tg} and the fifth equality follows since $T_g$ is an autoequivalence of $\bK$.
    
    Since $\Spc \bK^c$ stratifies $\bK$, that is, if $V^{BF}(A) \subseteq \Phi^{BF}(\bL)$ then $A \in \bL$ for any localizing ideal $\bL$ of $\bK$ and any $A \in \bK$, it now follows that $\bL$ is $G$-invariant if and only if $\Phi^{BF}(\bL)$ is $G$-invariant.
\end{proof}

When $\bK^c$ is stratified and centrally generated, one can provide a neat description of the homological spectrum of $\bK^c \rtimes G$. 

\bth{crossed-nos}
Let $\Spc \bK^c$ is Noetherian and $\bK$ be stratified by $\Spc \bK^c$. Moveover, 
assume that all Serre primes of $\Spch(\bK^c)$ are maximal (e.g., $\bK$ satisfies the Central Generation Hypothesis). Let $G$ be a group acting on $\bK$ by monoidal autoequivalences. Then $\bK \rtimes G$ satisfies the Nerves of Steel Condition.
\eth

\begin{proof}
    Let $G.\bP$ be an arbitrary point of $\Spc(\bK^c \rtimes G)$. The goal is to show that $G.\bP$ has a unique Serre ideal which is maximal with respect to the property that $\h^{-1}(\bS) = G.\bP$. Since we have assumed that $\bK$ is stratified by $\Spc(\bK^c)$ and satisfies the property that all Serre primes in $\Spch(\bK^c)$ are maximal, it follows from \thref{strat-nos} that $\bK$ satisfies Nerves of Steel; in particular, there is a unique homological prime $\bS$ which satisfies $\phi(\bS)=\bP$. 
    
    Let $E$ be the pure-injective which satisfies $\ker(\hat E \otimes -) = \overrightarrow{\bS},$ as in \thref{S-is-KerE}. Take 
    \[
    \bS' := \ker(\widehat{E \boxtimes e} \otimes -)^{\fp} \subseteq \modd(\bK^c \rtimes G).
    \]
    The embedding of $\bK^c$ as a direct summand in $\bK^c \rtimes G$ via $A \mapsto A \boxtimes e$ induces an embedding of $\modd(\bK^c)$ into $\modd(\bK^c \rtimes G),$ which is similarly denote by $F \mapsto F \boxtimes e.$ 
    
    If $F \in \bS$, then $F \boxtimes e \in \bS'$. Let $E'$ be the pure-injective of $\bK^c \rtimes G$ corresponding to $\bS'$. Then $E'$ has the form
    \[
    E' := \bigoplus_{g \in G} E'_g \boxtimes g
    \]
    for some collection of objects $E'_g \in \bK$. 

    By \coref{supp-pureinj}, the Balmer--Favi support of $E$ is $\{\bP\}$. Let $\bQ \not = \bP$ in $\Spc \bK^c$. We know that $\widehat{\Gamma_{\bQ} \unit} \in \overrightarrow{\bS}$; hence $\widehat{\Gamma_{\bQ} \unit} \boxtimes e \in \overrightarrow{\bS'}$. We now observe that
    \begin{align*}
    \widehat{\Gamma_{\bQ} \unit} \boxtimes e   \in \overrightarrow{\bS'} & \Rightarrow (\widehat{\Gamma_{\bQ} \unit} \boxtimes e) \otimes \hat E'    \cong 0\\
    & \Rightarrow (\widehat{\Gamma_{\bQ} \unit} \boxtimes e) \otimes (\widehat{E'_g \boxtimes g})    \cong 0 \quad \forall g \in G\\
    & \Rightarrow \Gamma_{\bQ} \unit \otimes E'_g \cong 0 \quad \forall g \in G\\
    & \Rightarrow \bQ \not \in V^{BF}(E_{g}') \quad \forall g \in G\\
    & \Rightarrow V^{BF}(E_g') = \{\bP\} \quad \forall g \in G\\
    & \Rightarrow V^{BF}(E')=\{G.\bP\},
    \end{align*}
    where the last two lines follow from \cite[Corollary 5.11]{CaiVashaw} and \cite[Proposition 3.8]{CaiVashaw}, respectively. 

    Note that if $\bS''$ was any Serre ideal of $\modd(\bK^c \rtimes G)$ maximal in the fiber of $G.\bP$, it would have a pure-injective $E''$ which would necessarily satisfy
    \[
    V^{BF}(E'') \supseteq V^{BF}(E'),
    \]
    by \thref{supp-pureinj-2}. There is at least one such Serre ideal, maximal in the fiber of $G.\bP$, such that the support of the pure injective is precisely $V^{BF}(E')$, since some ideal maximal in the fiber of $G.\bP$ must contain $\bS'$. Since $\bK \rtimes G$ is stratified by \thref{strat-crossprod}, it follows from \coref{max-loc-contain} that there is in fact a unique Serre ideal maximal in the fiber of $G.\bP$, and Nerves of Steel holds for $\bK^c \rtimes G$. 
\end{proof}

\bre{central-gen-g-fin}
If $\bK^c$ satisfies the Central Generation Hypothesis and $G$ is finite, then $\bK^c \rtimes G$ also satisfies central generation. Indeed, every thick ideal of $\bK^c \rtimes G$ is of the form
\[
\bI=\{ A \boxtimes g: A \in \bJ, g \in G\}
\]
for some thick ideal $\bJ$ of $\bK^c$. If $\{A_i : i \in I\}$ is a collection of central generators for $\bJ$, then 
\[
\left \{ \displaystyle{\bigoplus}_{g \in G} T_g(A_i) \boxtimes e : i \in I \right \}
\]
is a collection of central generators for $\bI$. 
\ere

Even when $\bK^c$ satisfies the Central Generation Hypothesis, if $G$ is infinite, $\bK^c \rtimes G$ might not satisfy Central Generation; we give an example illustrating this below. Hence, \thref{crossed-nos} produces the first examples where we can prove Nerves of Steel beyond Central Generation. 

\bex{central-gen-break}
Let $\bK$ denote the derived category of coherent sheaves on the complex affine plane $X:=\Spec \mathbb{C}[x,y]$. The compact part $\bK^c$ is the perfect derived category. The Balmer spectrum of $\bK^c$ is $X$, by Thomason's theorem \cite{Thomason}. Consider the group $G=\mathbb{Z}$ acting on the affine plane by horizontal translation; this induces an action of $G$ on $\bK^c$, and the induced action of $G$ on $\Spc \bK^c$ coincides with the original action of $G$ on $X$, by \cite[Example 3.5]{HuangVashaw}. 

The spectrum of $\bK^c \rtimes G$ was considered in \cite[Example 7.7]{HuangVashaw}, which can be used to show that $\bK^c \rtimes G$ does not satisfy the Central Generation Hypothesis. Thick ideals of $\bK^c \rtimes G$ are in bijection with $G$-orbits of closed subsets of $X$. Consider the thick ideal $\bI$ corresponding to the orbit of the origin. An indecomposable object will generate this ideal if and only if it is isomorphic to $M \boxtimes g$ for some $g \in G$, where $M$ is an object of $\bK^c$ whose support is a point $(a,0)$ for some fixed $a \in \mathbb{Z}$. But this object $M$ does not commute with any $\unit \boxtimes h$, as long as $h \not = 0$, since $T_h(M)$ has support $(a+h, 0)$, which is not the support of $M$. That is, no indecomposable generator of $\bI$ commutes with all objects of $\bK^c \rtimes G$; this is easily extended to non-indecomposable generators of $\bI$. 
\eex

\section{Finite Tensor Categories: Tale of Three Spectra} 
\label{9}
Recall that a {\em{finite tensor category}} \cite{Etingof-Ostrik} is an abelian rigid monoidal category $\bT$ over a field $k$ in which the endomorphism algebra of the unit object is reduced to scalars $\End(\unit)=k$ and such that $\bT$ is finite in the sense that it is equivalent to the category of finite dimensional modules $\modd{A}$ of a finite dimensional algebra $A$. Such categories play a prominent role in many areas of mathematics and mathematical physics. Every finite tensor category is Frobenius \cite[Section 2.8]{Etingof-Ostrik}, and we can form its stable category  
 $\underline{\bf T}$, which is a monoidal triangulated category. All such categories are the compact parts of rigidly-compactly generated monoidal triangulated categories (recall \deref{main-def}) and our treatment applies to them. This is proved in Appendix A of \cite{NVY3}, where it is shown that the indization (or ind completion) $\Ind(\bT)$ is a monoidal Frobenius abelian category, and its stable category $\underline{\Ind(\bT)}$ is a rigidly-compactly generated monoidal triangulated category with compact part $\underline{\bT}$. 

Important examples of finite tensor categories include the categories of finite dimensional modules $\bT= \modd(H)$ of finite dimensional Hopf algebras $H$. In this case, the above constructions work as follows. The indization of such a category $\Ind(\bT)$ is monoidally equivalent to the category of infinite dimensional modules $\Mod(H)$. Both $\modd(H)$ and $\Mod(H)$ are Frobenius abelian monoidal categories and the stable module category 
$\underline{\Mod}(H) \cong \underline{\bT}$ is 
a rigidly-compactly generated monoidal triangulated category with compact part $\underline{\modd}(H)$. 

 The results in this paper, along with the work in \cite{NVY3}, indicate that for $\underline{\bf T}$ there are three prominent spectra that have been introduced and developed by the authors in the setting of monoidal triangulated categories: 
\begin{itemize} 
\item[(i)] the homological spectra (cf. \deref{hom-spec}); 
\item[(ii)] the Balmer spectra \cite[Definition 3.2.1]{NVY1};
\item[(iii)] the cohomological spectra (via the categorical center)  \cite[Section 1.4]{NVY3}.
\end{itemize}
These spectra are important topological spaces that govern the geometry 
behind finite tensor categories. The reader should note that well-explored examples of finite tensor categories include representations for finite-dimensional Hopf algebras, including the small quantum groups, and finite group schemes. 

For finite group schemes\footnote{equivalently representations of finite-dimensional cocommutative Hopf algebra}, the study of these spectra 
started over 50 years with the pioneering work of Alperin \cite{Alperin}, Evens \cite{Evens}, Carlson \cite{Carlson, Carlson2}, Avrunin--Scott \cite{Avrunin-Scott}, 
Friedlander--Suslin \cite{Friedlander-Suslin}, and Friedlander--Pevtsova \cite{FriedlanderPevtsova} who proved important results involving the finite generation of the cohomology ring, and developed the theory of support varieties in connection 
with rank varieties ($\pi$-points) for finite groups, and more generally finite group schemes.

Let $\phi$ be the homological comparison map defined in Section~\ref{S:supportmap}. 
Furthermore, let $R^{\bullet}$ be the cohomology ring of $\bT$ and $C^{\bullet}$ be its categorical center (cf. \cite[Section 1.4]{NVY3}), and $\rho$ and $\eta$ be the cohomological support (continuous) comparison maps constructed in \cite{NVY3}. One has the following diagram of continuous maps that connect (i)-(iii): 
$$
\Spch(\underline{\bf T})\mathrel{\mathop{\rightarrow}^{\phi}}  \Spc(\underline{\bf T}) \mathrel{\mathop{\rightleftarrows}^{\rho}_\eta} \Proj(C^{\bullet}).
$$

\noindent
With the noncommutative version of Nerves of Steel (cf. 
\deref{nerves}), the authors propose the following capstone conjecture (updated from \cite[Conjecture E]{NVY3}) which give a clear picture of the interrelationships between these three spectra for finite tensor categories. 

\bcj{NVYconj} For every finite tensor category $\bT$, the continuous maps 
$\phi$, $\rho$ and $\eta$ 
$$
\Spch(\underline{\bT})\mathrel{\mathop{\rightarrow}^{\phi}}  \Spc(\underline{\bT}) \mathrel{\mathop{\rightleftarrows}^{\rho}_\eta} \Proj(C^{\bullet}).
$$
are homeomorphisms.
\ecj

In \cite{NVY3}, $\rho$ was defined as a continuous map 
$\rho: \Spc \underline{\bT} \to \Spech C^\bullet_{\underline{\bT}}$. Part of the conjecture 
is that its range is $\Proj(C^{\bullet})$. 
Furthermore, in \cite{NVY3} the continuous map $\eta : \Proj(C^{\bullet}) \to \Spc(\underline{\bT})$ was defined under the assumption that the central cohomological support map for $\underline{\bT}$ is a weak support datum in the sence of \cite{NVY1}. Part of the conjecture is that this holds for all finite tensor categories $\bT$.  

The conjecture has been verified for several important families of finite tensor categories using deep results in the literature, and the authors believe that Nerves of Steel should be true in this context. 

\bex{GroupScheme} {\bf [Finite Group Schemes]} Let $G$ be a finite group scheme over a field $k$. The module category for $G$ is equivalent to modules for a finite-dimensional cocommutative Hopf algebra $A$. Set $\bT=\modd(G)$, equivalently, $\modd(A)$. 
The cohomology ring $R^\bullet=\text{H}^{\bullet}(G,k)$ is finitely generated by a fundamental result of Friedlander-Suslin \cite{Friedlander-Suslin}. Moreover, the Etingof--Ostrik Conjecture holds in this case and $R^\bullet=C^\bullet$ since $\bT$ is symmetric. Using finite generation, the fact that $\rho$ and $\eta$ are homeomorphisms follows by work in \cite{Friedlander-Suslin} and by Benson--Iyengar--Krause--Pevtsova \cite{BIKP2018}. 

Finally, since $\bT$ is braided and rigid, the map $\phi$ is surjective. Stratification holds for finite group schemes by major results in \cite{BIKP2018}, and 
this can be used to show that $\phi$ is injective using pure injective modules (cf. \cite[Example 5.13]{Balmer2020c}). 
\eex
\vskip .15cm 
\bex{LieAlg} {\bf [Restricted Lie Algebras]} Let $({\mathfrak g},[p])$ be a finite-dimensional restricted Lie algebra over an algebraically closed field of characteristic $p>0$. Let $A=u({\mathfrak g})$ be the restricted enveloping algebra which is a finite-dimensional cocommutative Hopf algebra. Therefore, all the results from the prior example apply to this situation. 

More specifically, it is well-known that $\Proj(C^{\bullet})$ identifies with $\Proj(k[{\mathcal N}_{p}])$ where $k[{\mathcal N}_{p}]$ is the coordinate ring of the set of 
$p$-nilpotent elements ${\mathcal N}_{p}$. There are a natural set of tensor functors $\text{res}_{x}:\Mod(A)\rightarrow \Mod(u(\langle x \rangle)$ for every $x\in {\mathcal N}_{p}$ (after extending scalars). Each tensor functor gives rise to a homological prime 
${\mathcal B}_{x}$ for $x\in {\mathcal N}_{p}$. 

This is the most natural example of how the general theory works. However, the reader should be aware that finding such explicit tensor functors that detect nilpotence is a difficult problem. 

\eex
\vskip .15cm 
\bex{PW-Hopf} {\bf [Plavnik--Witherspoon Hopf Algebras]} \cjref{NVYconj} holds for the stable module category for an important class of finite-dimensional non-cocommutative Hopf algebras. These algebras were studied by Plavnik and Witherspoon \cite[Section 2]{Plavnik-Witherspoon}.  

Let $L$ be a finite group and $B$ be a finite-dimensional Hopf algebra with $H$ acting on $B$ by Hopf automorphisms. Set $A=(B^{*}\# kL)^{*}.$ As an algebra, $A \cong B \otimes k[L]$. 

The cohomology ring $R_A^{\bullet}$ of $A$ is finitely generated, but the support theory using $R_A^{\bullet}$ does not satisfy the tensor product property. (We use subscript $A$ here to emphasize that we work with the Hopf algebra $A$ and not with the Hopf algebra $B$.)

The categorical center, $C^{\bullet}$, is given by the following formula \cite[Theorem 5.2.1]{NVY3}: $C_{A}^{\bullet}\cong (C^{\bullet}_{B})^{L}$. The Balmer spectrum has been computed for the stable module category of finite-dimensional modules for $A$ (cf. \cite[Theorem 9.1.1]{NVY3}), and $\rho$ and $\eta$ are homeomorphisms under suitable conditions on $B$.  

There are two special cases of interest where $\rho$ and $\eta$ are homeomorphisms (cf. \cite[Propositions 9.2.1 and 9.3.1]{NVY3}). 

\begin{itemize}
\item[$\bullet$] Benson and Witherspoon studied a special case of this algebra in \cite{Benson-Witherspoon}. Let $G$ and $L$ be finite groups with $L$ acting on $G$ by group automorphisms. Set $A=(k[G] \# k L)^{*}$, where $k[G]$ is the dual of the group algebra of $G$, and $k L$ is the group algebra of $L$ where $k$ be a field of positive characteristic dividing the order of $G$. 

\item Let $\Omega$ be a finite group scheme, and let $k[\Omega]$ be the coordinate algebra of $\Omega$. Assume that $k$ is algebraically closed of characteristic $p$. 
One has $\Omega\cong \pi \ltimes \Omega_{0}$ where $\pi$ is a finite group and $\Omega_{0}$ is an infinitesimal group scheme. The Hopf algebra $k[\Omega]$ is a Plavnik--Witherspoon co-smash product Hopf algebra with $B=k[\Omega_{0}]$ and $L = \pi$. The monoidal triangular geometry associated to these algebras was studied by Negron--Pevtsova 
(cf.~ \cite[Section 10]{NegronPevtsovaSpec2023}).
\end{itemize} 

\cjref{NVYconj} holds for the stable module category for the Benson--Witherspoon Hopf algebras because $B=kG$ is a cocommutative Hopf algebra and the conditions of \thref{crossed-nos} hold.

In order to apply \thref{crossed-nos} for $k[\Omega]$, one has an algebra isomorphism 
with the group algebra of a direct product of cyclic groups of prime power order \cite[Section 9.3]{NVY3}. With this identification, one can use the argument for the group algebra of a finite group to prove stratification for the big stable module category for $k[\Omega]$ by an argument outlined below. The verification of 
the Central Generation Hypothesis has been done for several cases, see 
Negron--Pevtsova \cite[Section 10.E]{NegronPevtsovaSpec2023}. We do not know a general proof of central generation for $k[\Omega_0]$, although Negron--Pevtsova conjecture that it should hold \cite[Question 10.17]{NegronPevtsovaSpec2023}.

Nevertheless, we can deduce Nerves of Steel for $k[\Omega_0]$ and then for $k[\Omega]$ by the following argument, independent of central generation. Recall that $k[\Omega_0] \cong k G$ as algebras (although not necessarily as Hopf algebras) for $G$ an elementary abelian $p$-group. Since there is only one simple module, the tensor unit, in $\modd( k[\Omega_0])$, it follows that every thick subcategory of $\underline{\modd}(k[\Omega_0])$ is a thick ideal, and every localizing subcategory of $\underline{\Mod}(k[\Omega_0])$ is a localizing ideal, see \cite[Lemma 2.9, Lemma 3.13]{Stevenson2018}. It follows that, even though the tensor products on $\modd(kG)$ and $\modd(k [\Omega_0])$ are different, $\Spc \underline{\modd}(kG) \cong \Spc \underline{\modd}(k[\Omega_0])$, and the Balmer--Favi support is independent of the tensor product used. Since stratification of $\underline{\Mod}(kG)$ was achieved in \cite{BIK2011}, the big stable category of $k[\Omega_0]$ is also stratified. It remains to be shown that $\underline{\Mod}(k[\Omega_0])$ satisfies the Uniformity Hypothesis. 

Since the definition of pure-injective objects does not depend on a tensor product, the pure-injectives for $\underline{\Mod}(kG)$ coincide with the pure-injectives for $\underline{\Mod}(k[\Omega_0])$. Since $\underline{\Mod}(kG)$ is braided, all primes in $\Spch(\underline{\Mod}(kG))$ are maximal (\thref{thick-max}), and satisfy the property that if $\bS \in \Spch(\underline{\Mod}(kG))$ then $V^{BF}(E_{\bS}) = \{ \phi(\bS)\}$ (\coref{supp-pureinj}). 

Now let $\bP \in \Spc \underline{\modd}(k[\Omega_0])$. By the above discussion, we know that there is a pure-injective $E$ in $\underline{\Mod}(k[\Omega_0])$ such that $V^{BF}(E)=\{\bP\}$. We write $\otimes_1$ for the tensor product on $\underline{\Mod}(k[\Omega_0]) \cong \underline{\Mod}(kG)$ defined using the $kG$ coproduct, and write $\otimes_2$ for the tensor product defined using the $k[\Omega_0]$ coproduct. Consider 
\[
\bS:=\langle \ker(\hat E \otimes_1 -)^{\fp} \rangle_2,
\]
the smallest Serre ideal (with the ideal property taken with respect to the tensor product $\otimes_2$) of $\moddash{\underline{\modd}(k[\Omega_0])}$ containing $\ker(\hat E \otimes_1 -)^{\fp}.$ Since $\bS$ is a Serre ideal of $\moddash{\underline{\modd}(k[\Omega_0])}$, it has an associated pure-injective $E'$. Suppose $\bQ \not = \bP$. Since $\Gamma_{\bQ} \unit$ is in 
\[
\ker(\hat E \otimes_1 -) = \overrightarrow{\ker (\hat E \otimes_1 -)^{\fp}},
\]
and since $\overrightarrow{(-)}$ only depends on the triangulated structure of $\bK$ and not the tensor product, we have $\Gamma_{\bQ} \unit \in \overrightarrow{\bS}$, and so $\bQ \not \in V^{BF}(E')$. We claim that $\bS$ does not contain $\hat \unit$; this will imply that $E'$ cannot be isomorphic to 0, hence must have $\bP \in V^{BF}(E')$, by stratification. In fact, we argue that 
\[
\bS = \ker(\hat E \otimes_1-)^{\fp},
\]
from which it does indeed follow that $\hat \unit \not \in \bS$. Suppose $F \in \ker(\hat E \otimes_1 -)^{\fp}$, that is, $F$ is in $\moddash{\underline{\modd}(k[\Omega_0])}$ and $\hat E \otimes_1 F \cong 0$. We must argue that for any $G \in \moddash{\underline{\modd}(k[\Omega_0])}$, we have $\hat E \otimes_1 (F \otimes_2 G) \cong 0 \cong \hat E \otimes_1 (G \otimes_2 F);$ using the fact that $G$ is finitely-presented, it is enough to show that $\hat E \otimes_1 (F \otimes_2 \hat A) \cong 0 \cong \hat E \otimes_1 (\hat A \otimes_2 F)$ for all $A \in \underline{\modd}(k[\Omega_0]).$ This follows immediately from the fact that $\operatorname{Thick}(\unit)=\underline{\modd}(k[\Omega_0])$ and the fact that $\hat E \otimes_1 (F \otimes_2 \hat \unit) \cong \hat E \otimes_1 F \cong 0 \cong \hat E \otimes_1( \hat \unit \otimes_2 F)$.

To summarize, we have shown that $V^{BF}(E')=\{\bP\}$. Without loss of generality, we can assume that $\bS$ is in $\Spch(\bK^c)$ (potentially replacing $\bS$ by a Serre prime which is maximal in the fiber of $\bP$ containing it). 

By stratification, given any other pure-injective $E''$, if $\bP \in V^{BF}(E'')$, we would have $E' \in \Loc(E'')$, which would then imply that $\ker(\hat E' \otimes -) \supseteq \ker(\hat E'' \otimes -)$ by \coref{max-loc-contain}. Therefore, by $\ker(\hat E' \otimes -)$ is the unique maximal element in the fiber of $\bP$. That is, all homological primes are maximal, and $\underline{\Mod}(k[\Omega_0])$ satisfies the Uniformity Hypothesis, and hence Nerves of Steel for $\underline{\modd}(k[\Omega_0])$ follows, cf.~ \thref{strat-nos}. Since we have verified all necessary conditions for \thref{crossed-nos}, Nerves of Steel now also holds for $\underline{\modd}(k[\Omega])$. 
\eex

In the following theorem, we gather the major ideas in finite tensor categories that will be needed to prove \cjref{NVYconj}. 
The theorem is a combination of \cite[Theorem 8.2.1]{NVY3} and \thref{strat-nos}. 
One of the most well-known conjectures in the subject, is the Etingof--Ostrik Conjecture \cite[Section 1.5]{NVY3} which predicts that the cohomology ring 
for a finite tensor category ${\bf T}$, and its action of the cohomology between two objects for ${\bf T}$ is finitely generated. Note that 
conditions (i) and (ii) imply the Noetherianity of $\Spc(\underline{\bf T})$.

\bth{P:homprimes2} Let $\bT$ be a finite tensor category, ${\bT}^{\prime}$ be an ind-completion of $\bT$ and $\underline{\bT}^{\prime}$ be the corresponding stable module category. Assume that 
\begin{itemize}
\item[(i)] the Etingof--Ostrik Conjecture holds for the Drinfeld center of $\bT$; 
\item[(ii)] the central cohomological support of $\bf T$ has an extension to a faithful extended weak support datum to $\underline{\bT}^{\prime}$;
\item[(iii)] The Uniformity Hypothesis for $\underline{\bT}^{\prime}$ holds; 
\item[(iv)] $\underline{\bT}^{\prime}$ is stratified by $\Spc(\underline{\bf T}).
$
\end{itemize} 
Then the continuous maps 
$\phi$, $\rho$ and $\eta$ 
$$
\Spch(\underline{\bf T})\mathrel{\mathop{\rightarrow}^{\phi}}  \Spc(\underline{\bf T}) \mathrel{\mathop{\rightleftarrows}^{\rho}_\eta} \Proj(C^{\bullet})
$$
are homeomorphisms.
\eth

One of the long-standing problems has been to find rank variety analogs 
for finite-dimensional Hopf algebras (see \cite{pevtsova2009varieties} and \cite{NegronPevtsovaSpec2023}). The verification of \cjref{NVYconj} will show that there is a space of natural tensor functors arising from homological primes, like the situation for restricted Lie algebras, which governs the Balmer spectrum and the categorical center. See Balmer's remarks \cite[5.14 Remark]{Balmer2020c}. The question of the realization of these natural tensor functors, even in the cocommutative case, remains a mystery.  

\bibliography{bibl}
\bibliographystyle{myamsalpha}

\iffalse

\fi
\end{document}